\documentclass{amsart}
\usepackage{amssymb,a4wide}
\usepackage{graphicx} 
\usepackage{float} 
\usepackage{longtable} 
\usepackage{xcolor} 
\usepackage{comment}
\usepackage{eqparbox}
\newcommand{\Lb}[1]{\eqmakebox[Lhs][r]{$#1$}}
\newcommand{\Rb}[1]{\eqmakebox[Rhs][l]{$#1$}}

\def\pint{\langle \cdotp,\cdotp \rangle }
\def\R{\mathbb R}
\def\C{\mathbb C}
\def\N{\mathbb N}
\renewcommand{\b}{\mathfrak{b}}
\newcommand{\n}{\mathfrak{n}}
\renewcommand{\k}{\mathfrak{k}}
\newcommand{\g}{\mathfrak{g}}
\renewcommand{\a}{\mathfrak{a}}
\newcommand{\z}{\mathfrak{z}}
\newcommand{\h}{\mathfrak{h}}
\newcommand{\fu}{\mathfrak{u}}
\renewcommand{\v}{\mathfrak{v}}
\newcommand{\im}{\,\mathrm{Im}}
\newcommand{\Id}{\,\mathrm{Id}}
\DeclareMathOperator{\Span}{span}

\newtheorem{prop}{Proposition}[section]
\newtheorem{thm}[prop]{Theorem}
\newtheorem{lem}[prop]{Lemma}
\newtheorem{cor}[prop]{Corollary}

\theoremstyle{definition}
\newtheorem{defi}[prop]{Definition}
\newtheorem{rem}[prop]{Remark}
\newtheorem{exa}[prop]{Example}

\newtheorem{thmintro}{Theorem}

\title{(Non-)existence of complex structures on $2$-step nilpotent Lie algebras}
\author[M.L. Barberis, A. Moroianu, A. Tolcachier]{María Laura Barberis, Andrei Moroianu, Alejandro Tolcachier}

\address{(María L. Barberis and A. Tolcachier) 
 FAMAF, Universidad Nacional de C\'ordoba and CIEM, CONICET, 
  Av. Medina Allende s/n, Ciudad Universitaria \\ 
 X5000HUA C\'ordoba, Argentina.}

\email{mlbarberis@unc.edu.ar}
\email{atolcachier@unc.edu.ar}

\address{(A. Moroianu) Université Paris-Saclay, CNRS, Laboratoire de mathématiques d’Orsay, 91405, Orsay, France, and Institute of Mathematics “Simion Stoilow” of the Romanian Academy, 21 Calea Grivitei, 010702 Bucharest, Romania.}
\email{andrei.moroianu@math.cnrs.fr}

\begin{document}
\begin{abstract}
In this paper we determine which $8$-dimensional $2$-step nilpotent real Lie algebras admit (integrable) complex structures. We  construct explicit examples on $40$ out of the $43$ isomorphism classes, and show that there are obstructions to the existence of complex structures on the remaining three isomorphism classes. These obstructions are related to the set of elements of rank $2$ in the pencil defining the Lie algebra structure. We also show that for any even-dimensional Lie algebra $\h$, the direct sum $\g:=\R^2\oplus \h$ (and also some other, more general, extensions of $\R^2$ by $\h$) admits a complex structure  if and only if $\h$ admits a complex structure. In particular, this allows us to construct the first examples of $2$-step nilpotent Lie algebras in every even dimension $2n\geq8$ which do not admit complex structures. Finally, we show that every even-dimensional $2$-step nilpotent Lie algebra with $2$-dimensional center admits a complex structure.
\medskip


\medskip

\end{abstract}
\maketitle

\section{Introduction}

Nilmanifolds endowed with invariant complex structures are one of the main sources of examples of compact complex manifolds without K\"ahler metrics. Indeed, a nilmanifold $\Gamma\backslash N$, where $N$ is a simply connected nilpotent Lie group and $\Gamma$ is a lattice, admits a K\"ahler metric only when it is a torus \cite{BG}. For this reason, the existence of compatible Hermitian metrics with special properties on non-toral complex nilmanifolds has been extensively studied (see, for instance, \cite{AN,AV,COUV,EFV1,EFV2,FPS,FV,LoG,ST,UV}).
Since the complex structures under consideration are induced by left-invariant ones on $N$, the study of their existence is a purely algebraic problem on the Lie algebra $\n$ of $N$. Recall that a \emph{complex structure} on $\n$ is an endomorphism $J$ with $J^2=-\Id$ whose Nijenhuis tensor vanishes.

In dimension $4$, there are three isomorphism classes of nilpotent Lie algebras, only two of which admit complex structures, namely the abelian Lie algebra and $\R\oplus\h_3$, giving rise to the complex torus and the Kodaira--Thurston manifold, respectively. In his foundational work \cite{Sal}, Salamon carried out a detailed study of  nilpotent Lie algebras admitting complex structures and showed that, in dimension $6$, exactly $18$ out of the $34$ isomorphism classes of nilpotent Lie algebras admit complex structures. The concept of nilpotent complex structure, introduced in \cite{CFGU}, has attracted the attention of many authors since then (see, among others, \cite{GZZ, LU, LUV1, LUV2, LUV3, Rol,Zhang}). In dimension $8$, the nilpotent Lie algebras admitting non-nilpotent complex structures have been classified \cite{LU, LUV3}, but a complete classification of those admitting complex structures is not available. Moreover, in dimensions higher than $8$, no general criterion for the (non-)existence of complex structures is known.

In this paper, we study the existence of complex structures on $2$-step nilpotent Lie algebras. These Lie algebras appear in many different areas of mathematics. For instance, $2$-step nilpotent Lie groups are the simplest examples of Carnot groups, an important notion in analysis \cite{Fol}, and their left-invariant Riemannian geometry has been extensively studied \cite{Eb,Eb2}. A distinguished class among them is that of groups of type $H$, introduced by Kaplan \cite{Kap} as a generalization of the Heisenberg groups, whose solvable extensions, the so-called Damek--Ricci spaces, provide the first examples of non-symmetric harmonic Riemannian manifolds \cite{DR}. Moreover, naturally reductive homogeneous nilmanifolds are at most $2$-step nilpotent \cite{Gor}, and they were described in terms of representations of compact Lie groups in \cite{Lauret}. Some $2$-step nilpotent Lie algebras can also be naturally associated with graphs \cite{DM}, and complex structures on them were recently studied in \cite{AVe}. More generally, $2$-step nilpotent Lie groups have proved to be a fertile ground for the study of special geometric structures (see, among many others, \cite{ABB,AD,CdBM,dBM,dBMR,FGR,LRV,MO,MRV,OS}).

They also appear naturally in complex geometry: non-toral complex nilmanifolds admitting SKT metrics are necessarily $2$-step \cite{AN}, and so are, in real dimension~$8$, nilmanifolds carrying invariant hypercomplex structures \cite{DF} or (not necessarily invariant) astheno-K\"ahler metrics \cite{LU2}. Moreover, nilpotent Lie algebras of real dimension $2n$ admitting abelian complex structures with $h^{1,0}=n-1$ are $2$-step nilpotent \cite{PV}.

Complex structures on $2$-step nilpotent Lie algebras have been studied under additional assumptions, for instance when they are abelian \cite{ABD,AV} or bi-invariant \cite{De,DLV}. Recently, Barberis \cite{Bar} gave a structural characterization of $2$-step nilpotent Lie algebras admitting complex structures, by treating separately the cases where the complex structure is $2$-step or $3$-step nilpotent in the sense of \cite{CFGU}. However, deciding from this characterization whether a given $2$-step nilpotent Lie algebra admits a complex structure is not straightforward. Every even-dimensional $2$-step nilpotent Lie algebra of dimension at most $6$ admits a complex structure \cite{Bar,Sal}, but the question remained open in higher dimensions. We show that this fails already in dimension $8$, and, using the following result, which holds for arbitrary Lie algebras, we obtain $2$-step nilpotent Lie algebras without complex structures in every even dimension $\geq 8$.

\begin{thmintro}\label{thmA}
Let $\h$ be an even-dimensional Lie algebra with center $\z(\h)$ and let $\a$ be a $2$-dimensional vector space. For every element $\varphi\in\Lambda^2\a^*\otimes\z(\h)$, consider the Lie algebra $\g:=\a \oplus_\varphi \h$ whose underlying vector space is $\a \oplus \h$, with the following brackets:
\begin{equation*}
\begin{cases}
[e,f]_\g=\varphi(e,f), \qquad&\forall e,f\in\a,\\
[e,x]_\g=0,\qquad &\forall e\in\a, \forall x\in \h,\\
[x,y]_\g=[x,y]_\h, \qquad&\forall x,y\in\h.
\end{cases}
\end{equation*}
Then $\g$ admits a complex structure if and only if $\h$ admits a complex structure.
\end{thmintro}

The classification of $8$-dimensional $2$-step nilpotent real Lie algebras without abelian factors was recently obtained by Borovoi, Dina, and de Graaf \cite{BDG} (see also \cite{BR} for a classification over algebraically closed fields of characteristic different from $2$ and $3$). Using \cite{BDG} we prove the following:

\begin{thmintro}\label{thmB}
An $8$-dimensional $2$-step nilpotent Lie algebra admits a complex structure if and only if it is not isomorphic to one of
\begin{align*}
\n^5_{8,3}&=(0,0,0,0,0,e^{13},e^{24},e^{15}+e^{25}),\\
\n^6_{8,3}&=(0,0,0,0,0,e^{24},e^{15}+e^{23},e^{13}),\\
\n^7_{8,3}&=(0,0,0,0,0,e^{15}+e^{24},e^{14}+e^{23},e^{13}).
\end{align*}
Moreover, in every even dimension $\geq 10$ there exist both decomposable and indecomposable $2$-step nilpotent Lie algebras which do not admit complex structures.
\end{thmintro} 

To the best of our knowledge, these are the first examples of $2$-step nilpotent Lie algebras without complex structures. The three Lie algebras above have $3$-dimensional center, and the non-existence of complex structures follows from obstructions that we obtain in terms of the decomposable $2$-forms in their pencils, that is, in the image of $\n^*$ under the Chevalley--Eilenberg differential.

Finally, we obtain the following result for $2$-step nilpotent Lie algebras with $2$-dimensional center.

\begin{thmintro}\label{thmC}
Every even-dimensional $2$-step nilpotent Lie algebra with $2$-dimensional center admits a complex structure. Moreover, a large subclass, including all non-singular ones in the sense of \cite{Eb}, admits abelian complex structures.
\end{thmintro}

The proof uses the description of these Lie algebras in terms of pencil invariants \cite{Gau,Lau-Osc,LT}. We also determine which $8$-dimensional $2$-step nilpotent Lie algebras admit abelian complex structures (Theorem \ref{thm:abelian-dim8}).

\medskip
Let us now briefly describe the organization of the paper. After some preliminaries on $2$-step nilpotent Lie algebras, complex structures and Hermitian linear algebra in Section \ref{sec:prelim}, we prove Theorem \ref{thmA} in Section \ref{sec:3}, first in the case where $\varphi=0$, so $\g=\R^2\oplus \h$ is a direct sum of ideals (Theorem \ref{thm:g+r2}), then in general (Theorem \ref{thm:g+r2'}). We then turn our attention in Section \ref{sec:4} to $2$-step nilpotent $8$-dimensional Lie algebras and, using the recent classification results of \cite{BDG}, we show that all but three of them admit complex structures. In Section \ref{sec:5} we obtain a series of obstructions, in terms of the incidence relations in the pencils, to the existence of complex structures on $2$-step nilpotent $8$-dimensional Lie algebras with $3$-dimensional center, which turn out to be sufficient in order to rule out in Section \ref{sec:6}  the three remaining cases, thus proving the first part of Theorem \ref{thmB}. 
The second part follows from Theorem \ref{thmA} together with Proposition \ref{prop34}, which shows that the indecomposability is preserved by the construction of Theorem \ref{thmA}. This allows us to construct the first examples in every even dimension at least $8$ of $2$-step nilpotent Lie algebras which do not admit complex structures. These examples are indecomposable in dimension $8$, and can be chosen to be either decomposable or indecomposable in higher dimensions. In Section \ref{sec:7} we use the classification of $2$-dimensional pencils \cite{LT} to prove Theorem \ref{thmC}. Finally, in the Appendix we summarize the classification of $8$-dimensional $2$-step nilpotent Lie algebras, which was recently completed in \cite{BDG}.

\medskip
\noindent{\sc Acknowledgments.} This work was realized while the second named author was visiting FAMAF in the framework of the MATH AmSud Program SGP 24-MATH-12. The first and third named authors were partially supported by CONICET and SECyT-UNC. The authors thank Jorge Lauret for several interesting discussions, Adrián Andrada and Anna Fino for useful comments, and Luis Ugarte for his interest in our work and for pointing out reference \cite{Par} to us.

\medskip
\section{Preliminaries}\label{sec:prelim}

\subsection{Basics on $2$-step nilpotent Lie algebras}
A Lie algebra $\n$ is called $2$-step nilpotent if $[\n,\n']=0$ and $\n'\neq 0$, where $\n'=[\n,\n]$ denotes its commutator ideal. Equivalently, $\n$ is 2-step nilpotent if $0\neq \n'\subset \z$, where $\z$ denotes the center of $\n$.

Let $\n$ be a $2$-step nilpotent Lie algebra and set $r:=\dim(\n')$. Then, $\n$ determines an $r$-dimensional subspace  $P_\n\subset \Lambda^2(\n/\n')^*$, called a {\em pencil}, which is equal to the image of $\n^*$ under the Chevalley--Eilenberg differential (note that $d(\n^*)$ is contained in the subspace $\Lambda^2(\n/\n')^*\subset \Lambda^2(\n)^*$ since $z\lrcorner d\eta=0$ for every $\eta\in\n^*$ and $z\in \n'$). 

The pencils of $2$-step nilpotent Lie algebras can be described more concretely in the presence of a scalar product.
Given a $2$-step nilpotent Lie algebra $\n$ with an inner product $ \pint$, consider the orthogonal decomposition $\n=\n' \oplus \v$. We recall that the Lie bracket on $\n$ is determined by a linear map $j: \n ' \to \mathfrak{so}(\v)$ defined by
\begin{equation}\label{eq:jz-maps}    
 \langle j(z) v, w\rangle := \langle z, [v,w] \rangle, \qquad z\in \n ', \, v,w\in \v . \end{equation}
 The notation $j (z)$ was first introduced
by A. Kaplan in \cite{Kap} to study $2$-step nilpotent Lie groups  of Heisenberg type.

\subsection{Complex structures on Lie algebras}
An {\em almost complex structure} on a Lie algebra $\g$ is an endomorphism $J$ of $\g$ satisfying $J^2=-\Id$. The Nijenhuis tensor of $J$ is defined by 
$$N_J(x,y):=[x,y]-[Jx,Jy]+J([Jx,y]+[x,Jy]),\qquad\forall x,y\in\g.$$
Clearly, 
\begin{equation}\label{nij}
    N_J(x,x)=N_J(x,Jx)=0,\qquad N_J(Jx,y)=-JN_J(x,y),\qquad \forall x,y\in\g.
\end{equation}
An almost complex structure $J$ is called {\em integrable} if $N_J=0$. An integrable almost complex structure is called a {\em complex structure}.

Every almost complex structure determines a decomposition $\g^*\otimes\C=\Lambda^{1,0}\g^*\oplus\Lambda^{0,1}\g^*$, with
\[
\Lambda^{1,0}\g^*:=\{\eta\in \g^*\otimes\C\ |\ \eta\circ J=i\eta\},\qquad \Lambda^{0,1}\g^*:=\{\eta\in \g^*\otimes\C\ |\ \eta\circ J=-i\eta\}.
\]
Then $J$ is integrable if and only if \[d(\Lambda^{1,0}\g^*)\subset\Lambda^{1,1}\g^*\oplus\Lambda^{2,0}\g^*:=\Lambda^{1,0}\g^*\wedge(\g^*\otimes\C),\]
where $d$ denotes the Chevalley--Eilenberg differential.

Recall that an almost complex structure $J$ on $\g$ is called \emph{abelian} if $[Jx,Jy]=[x,y]$ for all $x,y\in\g$. An abelian almost complex structure is automatically integrable and preserves the center. In particular, $\dim\z$ is even, and the orthogonal complement $\z^\perp$ with respect to any Hermitian metric is $J$-invariant.

\subsection{Hermitian $4$-dimensional linear algebra}\label{sec23}

Let $\v_0$ be a 4-dimensional vector space endowed with a Hermitian structure $(g,J)$. We introduce the direct sum decomposition $$\Lambda^2\v_0^*=L^+(\v_0)\oplus L^-(\v_0),$$ where
\begin{equation}\label{lambda}
    L^\pm(\v_0)=\{\alpha\in \Lambda^2\v_0^*\ |\ \alpha(x,y)=\pm\alpha(Jx,Jy),\ \forall x,y\in\v_0\}.
\end{equation}

The complex structure $J$ extends to the dual space $\v_0^*$ by $J\eta:=-\eta\circ J$, and then to the exterior algebra $\Lambda\v_0^*$ as an automorphism. By \eqref{lambda}, this extension acts as the identity on $L^+(\v_0)$ and minus the identity on $L^-(\v_0)$.

\begin{rem}\label{rem2}
    Identifying $\Lambda^2\v_0^*$ with skew-symmetric endomorphisms of $\v_0$, $L^+(\v_0)$ and $L^-(\v_0)$ correspond to endomorphisms commuting, respectively anti-commuting with $J$. In this way, $J$ acts by left composition on $L^-(\v_0)$ by $L^-(\v_0)\ni\beta\mapsto J\circ\beta\in L^-(\v_0)$, where 
    \begin{equation}\label{jbe}
    (J\circ\beta)(x,y):=-\beta(Jx,y).
    \end{equation}
\end{rem}

For later use, we prove here the following standard result:

\begin{lem}
    For every $\eta\in\v_0^*$ and $\beta\in L^-(\v_0)$, one has 
  \begin{equation}\label{eb}
      \eta\wedge\beta=J\eta\wedge (J\circ\beta),\qquad J\eta\wedge\beta+\eta\wedge (J\circ\beta)=0.
  \end{equation}
\end{lem}

\begin{proof}
    With respect to the orientation induced by $J$ on $\v_0$, elements of $L^-(\v_0)$ are self-dual 2-forms. 
    Moreover, for every $x,y\in\v_0$ we have by \eqref{jbe}
    $$(Jx\lrcorner (J\circ\beta))(y)=(J\circ\beta)(Jx,y)=\beta(x,y)=(x\lrcorner\beta)(y),$$
    showing that $Jx\lrcorner (J\circ\beta)=x\lrcorner\beta$. This, together with the standard formula 
    $$*(\eta\wedge\beta)=\eta^\sharp\lrcorner *\beta$$
    relating the Hodge duality with exterior products and contractions, yields
    $$*(J\eta\wedge (J\circ\beta))=J\eta^\sharp\lrcorner *(J\circ \beta)=J\eta^\sharp\lrcorner (J\circ \beta)=\eta^\sharp\lrcorner \beta=*(\eta\wedge\beta),$$
    thus proving the first part of \eqref{eb}. The second one follows by replacing $\eta$ with $J\eta$.
\end{proof}

\medskip
\section{Complex structures on certain extensions of Lie algebras}\label{sec:3}

We start by proving in Theorem \ref{thm:g+r2} a particular case of Theorem \ref{thmA} stated in the introduction. The general case (Theorem \ref{thm:g+r2'}) will be obtained with similar arguments, except for one point that requires some extra care. 

\begin{thm}\label{thm:g+r2}
Let $\h$ be an even-dimensional Lie algebra and let $\a$ be a $2$-dimensional abelian Lie algebra. Then the Lie algebra $\g:=\a \oplus \h$, the direct sum of ideals $\a$ and $\h$, admits a complex structure if and only if $\h$ admits a complex structure.
\end{thm}
\begin{proof}
    The converse direction is obvious. For the direct implication, let $J$ be a complex structure on $\g$. If $J\h=\h$, then the restriction of $J$ to $\h$ is clearly integrable.  If $J\a=\a $, one can define an almost complex structure on $\g/{\a }\simeq\h$ by $\bar J(x+\a ):=Jx+\a $, which is again clearly integrable.
    
    We thus assume from now on that $J\h\neq\h$ and $J\a \neq\a $. Then $\a \cap J\a =0$, and $\h_0:=\h\cap J\h$ has codimension 2 in $\h$ and is $J$-invariant. 

    Let $\{e_1,e_2\}$ be a basis of $\a $ and decompose 
    \begin{equation}\label{je}
        Je_1=f_1+\xi_1,\qquad Je_2=f_2+\xi_2,\qquad\hbox{with }f_i\in\a , \ \xi_i\in\h\ \hbox{for }i=1,2.
    \end{equation}
    Note that $\xi_1$ and $\xi_2$ are linearly independent. Indeed, if $\lambda\,\xi_1+\mu\,\xi_2=0,$ then $J(\lambda\, e_1+\mu\, e_2)\in\a $, so $\lambda \,e_1+\mu \,e_2\in \a \cap J\a =0$, whence $\lambda=\mu=0$.

   Since all brackets of $\g$ belong to $\h$, the integrability of $J$ shows that for every $x,y\in \g$
    \begin{equation}\label{xy}
        [x,y]-[Jx,Jy]\in\h_0,\qquad [Jx,y]+[x,Jy]\in\h_0,
    \end{equation}
    and for every $x\in \g$ and $i=1,2$, 
    $$0=N_J(e_i,x)=-[Je_i,Jx]+J[Je_i,x]=-[\xi_i,Jx]+J[\xi_i,x].$$
    We thus have for $i=1,2$:
    \begin{equation}\label{xi}
        J[\xi_i,x]=[\xi_i,Jx],\qquad\forall x\in\g,
    \end{equation}
    which also implies 
    \begin{equation}\label{xi1}
        [\xi_i,x]\in\h_0,\qquad\forall x\in\g.
    \end{equation}

We will now distinguish two cases. 
\medskip

{\bf Case 1.} Assume that $\Span(\xi_1,\xi_2)\cap\h_0=0$. Then $\h=\Span(\xi_1,\xi_2)\oplus\h_0$.
Consider the unique almost complex structure $J_\h$ on $\h$ satisfying
    $$ J_\h\xi_1=\xi_2,\qquad\hbox{and}\qquad J_\h x=Jx,\ \forall x\in\h_0.$$
    Using \eqref{xy}, we compute for every $x,y\in\h_0$
     \begin{equation}\label{case1}
    \begin{split}
        N_{J_\h}(x,y)&=[x,y]-[J_\h x,J_\h y]+ J_\h([ J_\h x,y]+[x, J_\h y])\\
        &=[x,y]-[Jx,Jy]+ J_\h([Jx,y]+[x,Jy])\\
        &=[x,y]-[Jx,Jy]+J([Jx,y]+[x,Jy])=N_J(x,y)=0.
        \end{split}
   \end{equation}
    Similarly, using \eqref{xi} and \eqref{xi1}, we compute for every $x\in\h_0$
    \begin{eqnarray*}
        N_{J_\h}(\xi_1,x)&=&[\xi_1,x]-[ J_\h \xi_1, J_\h x]+J_\h([ J_\h \xi_1,x]+[\xi_1, J_\h x])\\
        &=&[\xi_1,x]-[\xi_2,Jx]+ J_\h([\xi_2,x]+[\xi_1,Jx])\\
        &=&-J(J[\xi_1,x]-[\xi_1,Jx])+J[\xi_2,x]-[\xi_2,Jx]=0.
    \end{eqnarray*}
   By \eqref{nij} we thus get $N_{J_\h}=0$, so $J_\h$ is a complex structure on $\h$.

\medskip
{\bf Case 2. } Assume, finally, that $\Span(\xi_1,\xi_2)\cap\h_0\neq0$. We claim that in this case $\Span(\xi_1,\xi_2)\subset\h_0$. 
Consider some linear combination $a\,\xi_1+b\,\xi_2\in \h_0\setminus\{0\}$. By \eqref{je} we get 
\begin{equation}\label{ja}
    J(a\,e_1+b\,e_2)=f+\xi
\end{equation}
with $f\in\a $ and $\xi =a\,\xi_1+b\,\xi_2\in\h_0\setminus\{0\}$. Then 
\begin{equation}\label{jb}
    Jf=-(a\,e_1+b\,e_2)-J\xi. 
\end{equation}
Note that $f$ cannot be proportional to $a\,e_1+b\,e_2$: indeed, if $f=c\,(a\,e_1+b\,e_2)$, \eqref{ja} and \eqref{jb} imply that $c\,(f+\xi)=-(a\,e_1+b\,e_2)-J\xi$, whence $J\xi=-c\,\xi$. Since $J$ has no real eigenvalues, this contradicts the fact that $\xi\neq 0$. Thus $f$ and $a\,e_1+b\,e_2$ are linearly independent elements of $\a $, so \eqref{ja} and \eqref{jb} show that $J\a \subset \a \oplus \h_0$, proving our claim. In particular, $\a\oplus\h_0$ is $J$-invariant.

 Since $(\a\oplus\h_0,J)$ is a complex hyperplane in $(\g,J)$ viewed as complex vector space, every vector not contained in it generates a complementary complex line. Let $v_1\in \h\setminus \h_0\subset \g\setminus (\a\oplus\h_0)$ and set $\mathfrak b :=\Span (v_1, Jv_1)$, which is thus a $J$-invariant complement of $\a \oplus \h_0$ in $\g$.
Consider the projection $\pi_\h : \g\to \h$ along $\a$. Let $\h_1:=\pi_\h (\b) $, then $\h= \pi_\h (\g)=\h_0+\h_1$. For dimensional reasons, we have that $\pi_\h :\b \to \h_1$ is an isomorphism, therefore, setting $v_2:=\pi_\h (Jv_1)$, it follows that $\{v_1, v_2\}$ is a basis of $\h_1$.

Let us now define the almost complex structure $J_\h$ on $\h$ by 
    $$ J_\h v_1:=v_2,\qquad J_\h v_2:=-v_1,\qquad J_\h x:=Jx\qquad\forall x\in\h_0.$$
    The very same computation as in \eqref{case1} shows that $N_{J_\h}(x,y)=0$ for every $x,y\in\h_0$.
    
    On the other hand, set $(Jv_1)_\a:=Jv_1-\pi_\h (Jv_1) \in \a $, so that $J_\h v_1=Jv_1- (Jv_1)_\a$.   Using the second part of \eqref{xy} applied to $x=v_1$, which tells us that $[Jv_1,y]+[v_1,Jy]\in\h_0,$ we compute for every $y\in\h_0$:
    \begin{eqnarray*}
        N_{J_\h}(v_1,y)&=&[v_1,y]-[ J_\h v_1, J_\h y]+ J_\h([ J_\h v_1,y]+[v_1, J_\h y])\\
        &=&[v_1,y]-[Jv_1,Jy]+ J_\h([Jv_1,y]+[v_1,Jy])\\
        &=&N_{J}(v_1,y)=0.
    \end{eqnarray*}
   By \eqref{nij} we also get $N_{J_\h}( J_\h v_1,x)=0$ for every $x\in\h_0$. Since $\Span(v_1, J_\h v_1)\oplus\h_0=\h$, we obtain $N_{J_\h}=0$, so $J_\h$ is a complex structure on $\h$.
\end{proof}

\begin{exa}
    Let $\mathfrak f_4$ be the (unique up to isomorphism) $3$-step nilpotent $4$-dimensional Lie algebra. It is well known (cf. \cite[Proposition 2.1, Corollary 2.3]{Snow} and  \cite[Proposition 4.1]{GR}) that $\mathfrak f_4$ admits no complex structure. As a direct consequence of Theorem \ref{thm:g+r2} we thus obtain that $\R^{2k}\oplus\mathfrak f_4$ admits no complex structure for $k\geq 1$, a fact which was proved in \cite[Proposition 3.26]{Par}.
\end{exa}

Notice that the bracket of $\a$ was never used in the proof except in the case where $J\a=\a$. In fact, with a slight modification, we can prove the following more general result:

\begin{thm}\label{thm:g+r2'}
    Let $\h$ be an even-dimensional Lie algebra with center  $\z(\h)$ and let $\a$ be a $2$-dimensional vector space. For every element $\varphi\in\Lambda^2\a^*\otimes\z(\h)$, consider the Lie algebra $\g:=\a \oplus_\varphi \h$ whose underlying vector space is
    $\a \oplus \h$, with the following brackets:
    \begin{equation}
        \begin{cases}
            [e,f]_\g=\varphi(e,f),
            \qquad&\forall e,f\in\a,\\
            [e,x]_\g=0,\qquad &\forall e\in\a, \forall x\in \h,\\
            [x,y]_\g=[x,y]_\h,
            \qquad&\forall x,y\in\h.\\
        \end{cases}
    \end{equation}
    Then $\g$
    admits a complex structure if and only if $\h$ admits a complex structure.
\end{thm}
\begin{proof}
    Note that the condition that $\varphi$ takes values in the center of $\h$ is necessary and sufficient for $[\cdot,\cdot]_\g$ to satisfy the Jacobi identity.

    The converse implication is clear: if $J_\h$ is a complex structure on $\h$, then for every choice of complex structure $J_\a$ on the vector space $\a$, the almost complex structure $J:=J_\a\oplus J_\h$ on $\g$ is integrable. 

    The proof of the direct statement is exactly the same as the one of Theorem \ref{thm:g+r2}, except for the case where $J\a=\a$. Indeed, assume that $\varphi\neq0$ (if $\varphi=0$ we are in the framework of Theorem \ref{thm:g+r2}). Since $\a$ is no longer an ideal, one cannot identify $\h$ with $\g/\a$ anymore. However, if $J$ is an integrable complex structure on $\g$ with $J\a=\a$, then for every $e\in \a$ and $x\in\h$, denoting $Jx=f+y$ with $f\in\a$ and $y\in\h$, we have 
    \begin{eqnarray*}
        0&=&N_J(e,x)=[e,x]_\g-[Je,Jx]_\g+J([e,Jx]_\g+[Je,x]_\g)\\
        &=&-[Je,Jx]_\g+J[e,Jx]_\g=-\varphi(Je,f)+J\varphi(e,f).
    \end{eqnarray*}
     Since $\varphi\neq 0$, taking $e:=f$ in the above calculation shows that $\varphi(Jf,f)=0$, whence $f=0$. This holds for every $x\in \h$, so $J\h\subset\h$, which implies that $J_\h:=J|_\h$ is a complex structure on $\h$.
\end{proof}

In particular, consider an odd-dimensional Lie algebra $\k$. Applying repeatedly
Theorem \ref{thm:g+r2'} to the Lie algebra $\h:=\k\oplus\R\xi$, 
and to $\varphi:=e^{12}\otimes\xi\in\Lambda^2(\R^2)^*\otimes \R\xi$, we obtain the following:

\begin{cor}
Let $\k$ be an odd-dimensional Lie algebra and $n\geq 1$. Then $\k\oplus\R$ admits a complex
structure if and only if $\k\oplus\h_{2n+1}$ admits a complex structure.
\end{cor}

\smallskip

\begin{rem}    
The Lie algebra $\g$ constructed in Theorem \ref{thm:g+r2'} fits in an exact sequence 
$$0\to\h\to\g\to\R^2\to 0,$$
so it is an extension of the abelian Lie algebra $\R^2$ by $\h$. In \S\ref{s65} we will construct a complex structure on an extension $\n$ of $\R^2$ by an 8-dimensional Lie algebra $\n_0$ which does not admit complex structures, thus showing that Theorem \ref{thm:g+r2'} does not hold for arbitrary extensions of $\R^2$.
\end{rem}

We end this section with a result  about decomposability. Recall that a Lie algebra $\g$ is called decomposable if there exist two complementary non-zero ideals $\g_1$ and $\g_2$ in $\g$.
Our aim is to show the following:
\begin{prop}\label{prop34}
    Let $\h$ be an indecomposable Lie algebra with center  $\z(\h)$ and let $\a$ be a $2$-dimensional vector space. Then for every non-zero element $\varphi\in\Lambda^2\a^*\otimes\z(\h)$, the Lie algebra $\g:=\a \oplus_\varphi \h$ defined in Theorem \ref{thm:g+r2'} is indecomposable.
\end{prop}

\begin{proof}
    Assume that $\g=\g_1\oplus \g_2$ is a direct sum of non-trivial ideals. We proceed in two steps.

    \medskip
    {\bf Step 1.} We show that $\g_1$ and $\g_2$ are not contained in $\h$. Assume that $\g_1\subset \h$. We have $\dim(\g_2\cap\h)\geq\dim(\g_2)-2$, so $\dim(\g_1)+\dim(\g_2\cap\h)\geq\dim(\g)-2=\dim\h$. Since moreover 
    $\g_1\cap(\g_2\cap\h)=0$, we obtain $\h=\g_1\oplus(\g_2\cap\h)$. This is clearly a direct sum of ideals of $\h$. The indecomposability of $\h$ and the fact that $\g_1\neq0$ then implies $\g_1=\h$ and thus $\dim(\g_2)=2$.
    Let $\pi_\a$ denote the projection on $\a$. We have $\a=\pi_\a(\g_1)+\pi_\a(\g_2)=\pi_\a(\g_2)$, 
    and therefore $\g_2=\Span(e+x,f+y)$, where $\{e,f\}$ is a basis of $\a$ and $x,y\in\h$. In particular, we have $\varphi(e,f)\neq 0$. Now, since $[\g_2,\g_1]=0$, we have that $[e+x,\h]=[f+y,\h]=0$, whence $x,y\in \z(\h)$ (because $e$ and $f$ commute with $\h$). Consequently, $[e+x,f+y]=[e,f]=\varphi(e,f)\in \h\setminus\{0\}$, contradicting the fact that $[e+x,f+y]\in\g_2$. A similar argument shows that $\g_2$ is not contained in $\h$.

       \medskip 

    {\bf Step 2.} Let $\pi_\a$ denote the projection on $\a$. By Step 1 we have $\pi_\a(\g_i)\neq 0$ for $i=1,2$. Since moreover $\pi_\a(\g_1)+\pi_\a(\g_2)=\a$, there exist vectors $e+x\in\g_1$ and $f+y\in\g_2$ such that $\{e,f\}$ is a basis of $\a$ and $x,y\in\h$. Denote by $\xi:=[e,f]\neq 0$. We have $[e+x,f+y]\in[\g_1,\g_2]=0$, so $[x,y]=-\xi$. But then $\xi=-[x,y]=-[x,f+y]\in\g_2$, and $\xi=-[x,y]=-[e+x,y]\in\g_1$, which implies $\xi\in\g_1\cap\g_2=0$, contradicting the fact that $\xi\neq 0$.
\end{proof}

\medskip
\section{Construction of complex structures}\label{sec:4}

The classification of $8$-dimensional $2$-step nilpotent Lie algebras has been recently obtained by Borovoi, Dina and de Graaf \cite{BDG}. We refer to the Appendix for the notation of the different isomorphism classes.
 
Based on this classification, we prove in this section that every $8$-dimensional $2$-step nilpotent Lie algebra $\n$ admits a complex structure, except for the isomorphism classes $\n_{8,3}^5$, $\n_{8,3}^6$ and $\n_{8,3}^7$. This will prove the existence part of Theorem \ref{thmB}.

By the classification of Salamon \cite{Sal} (see also \cite[Section~4.2]{Bar}), every $6$-dimensional $2$-step nilpotent Lie algebra admits a complex structure. Taking the product with the standard complex structure on $\R^2$, we obtain complex structures on $\R^5\oplus\h_3$, $\R^3\oplus\h_5$, $\R^3\oplus\n_{5,2}$, $\R^2\oplus\h_3\oplus\h_3$, $\R^2\oplus\h_3^\C$, $\R^2\oplus\n_{6,2}$ and $\R^2\oplus\n_{6,3}$.

For each of the remaining Lie algebras, we give an explicit complex structure $J$ in Tables~\ref{tab:complex-1}, \ref{tab:complex-2}, \ref{tab:complex-3} and \ref{tab:complex-4}, organized according to $\dim(\n')$. We define $J$ on $1$-forms by $J\alpha:=-\alpha\circ J$, and we describe it by giving a basis of $\Lambda^{1,0}\n^*=\{\alpha+iJ\alpha \mid \alpha\in\n^*\}$.

Recall that $J$ is integrable if and only if $d(\Lambda^{1,0}\n^*)\subset\Lambda^{2,0}\n^*\oplus\Lambda^{1,1}\n^*$, which can be checked directly in each case. In the tables we only write the non-zero differentials.

\begin{exa}
     Consider, for instance, the Lie algebra $\n:=\n_{8,3}^2$ defined by 
$$de^6=e^{14}+e^{35},\qquad de^7=e^{15}+e^{24}, \qquad de^8=e^{23}.$$
Define the almost complex structure $J$ whose $(1,0)$ space $\Lambda^{1,0}\n^*$ is spanned by 
$$\sigma_1:=e^1+ie^4,\qquad \sigma_2:=e^2+ie^5,\qquad \sigma_3:=e^3+ie^7,\qquad \sigma_4:=e^6+ie^8.$$
We have 
$$d\sigma_1=d\sigma_2=0,\qquad d\sigma_3=i(e^{15}+e^{24})=-i \im(\sigma_1\wedge\bar\sigma_2)\in\Lambda^{1,1}\n^*,\text{ and}$$
$$d\sigma_4=e^{14}+e^{35}+ie^{23}=\frac{i}2\sigma_1\wedge\bar\sigma_1+\frac{i}2\sigma_2\wedge(\sigma_3+\bar\sigma_3)\in\Lambda^{1,1}\n^*\oplus\Lambda^{2,0}\n^*,$$
thus showing that $J$ is integrable.
\end{exa}
With similar arguments one can easily check the integrability of all complex structures described in the tables below.

{\small
\renewcommand{\arraystretch}{1.3}
\begin{longtable}{|l|l|l|}
\hline
$\n$ & $de^8$ & $\Lambda^{1,0}\n^*$ \\
\hline
\endfirsthead
\hline
$\n$ & $de^8$ & $\Lambda^{1,0}\n^*$ \\
\hline
\endhead
\multicolumn{3}{r}{\footnotesize\emph{Continued on next page}} \\
\endfoot
\caption{Complex structures on $8$-dimensional $2$-step nilpotent Lie algebras with $\dim(\n')=1$ not of the form $\R^2\oplus\mathfrak m$}
\label{tab:complex-1}
\endlastfoot
$\mathbb R \oplus\mathfrak h_7$ & $e^{12}+e^{34}+e^{56}$ & $\Span( e^1+ie^2,\, e^3+ie^4,\, e^5+ie^6,\, e^7+ie^8)$ \\
\hline
\end{longtable}
}

{\small
\renewcommand{\arraystretch}{1.3}
\begin{longtable}{|l|l|l|}
\hline
$\n$ & $(de^7,de^8)$ & $\Lambda^{1,0}\n^*$ \\
\hline
\endfirsthead
\hline
$\n$ & $(de^7,de^8)$ & $\Lambda^{1,0}\n^*$ \\
\hline
\endhead
\multicolumn{3}{r}{\footnotesize\emph{Continued on next page}} \\
\endfoot
\caption{Complex structures on $8$-dimensional $2$-step nilpotent Lie algebras with $\dim(\n')=2$ not of the form $\R^2\oplus\mathfrak m$}
\label{tab:complex-2}
\endlastfoot
$\mathbb R \oplus\mathfrak n_{7,2,A}$ & $(e^{12},\,e^{14}+e^{35})$ & $\Span( e^1+ie^6,\, e^2+ie^4,\, e^3+ie^5,\, e^7+ie^8)$ \\
\hline
$\mathbb R \oplus\mathfrak n_{7,2,B}$ & $(e^{12}+e^{34},\,e^{15}+e^{23})$ & $\Span( e^1+i(e^5-e^4),\, e^2+ie^4,\, e^3+ie^6,\, e^7+ie^8)$ \\
\hline
$\mathfrak h_3\oplus\mathfrak h_5$ & $(e^{12},\,e^{34}+e^{56})$ & $\Span( e^1+ie^2,\, e^3+ie^4,\, e^5+ie^6,\, e^7+ie^8)$ \\
\hline
$\mathfrak n_{8,2}^{1}$ & $(e^{12}+e^{56},\,e^{34}+e^{56})$ & $\Span( e^1+ie^2,\, e^3+ie^4,\, e^5+ie^6,\, e^7+ie^8)$ \\
\hline
$\mathfrak n_{8,2}^{1\text{-bis}}$ & $(e^{12}-e^{36}+e^{45},\,e^{34}-e^{56})$ & $\Span( e^1+ie^2,\, e^3+ie^4,\, e^5+ie^6,\, e^7+ie^8)$ \\
\hline
$\mathfrak n_{8,2}^{2}$ & $(e^{15}+e^{46},\,e^{14}+e^{23})$ & $\Span( e^1+ie^4,\, e^2+ie^3,\, e^5+ie^6,\, e^7+ie^8)$ \\
\hline
$\mathfrak n_{8,2}^{3}$ & $(e^{15}+e^{26}+e^{34},\,e^{14}+e^{23})$ & $\Span( e^1+ie^2,\, e^3-ie^4,\, e^5+ie^6,\, e^7+ie^8)$ \\
\hline
$\mathfrak n_{8,2}^{4}$ & $(e^{13}+e^{25},\,e^{16}+e^{24})$ & $\Span( e^1+ie^2,\, e^3+ie^5,\, e^4-ie^6,\, e^7+ie^8)$ \\
\hline
$\mathfrak n_{8,2}^{6}$ & $(e^{16}+e^{25}+e^{34},\,e^{12})$ & $\Span( e^1+ie^2,\, e^3+ie^4,\, e^5-ie^6,\, e^7+ie^8)$ \\
\hline
\end{longtable}
}

{\small
\renewcommand{\arraystretch}{1.3}
\begin{longtable}{|l|l|l|}
\hline
$\n$ & $(de^6,de^7,de^8)$ & $\Lambda^{1,0}\n^*$ \\
\hline
\endfirsthead
\hline
$\n$ & $(de^6,de^7,de^8)$ & $\Lambda^{1,0}\n^*$ \\
\hline
\endhead
\multicolumn{3}{r}{\footnotesize\emph{Continued on next page}} \\
\endfoot
\caption{Complex structures on $8$-dimensional $2$-step nilpotent Lie algebras with $\dim(\n')=3$ not of the form $\R^2\oplus\mathfrak m$}
\label{tab:complex-3}
\endlastfoot
$\mathbb R \oplus\mathfrak n_{7,3,A}$ & $(e^{12},\,e^{23},\,e^{24})$ & $\begin{array}[t]{@{}l@{}} \Span( e^1+ie^2,\, e^3+ie^4,\\ \phantom{\Span(}e^5+ie^6,\, e^7+ie^8) \end{array}$ \\
\hline
$\mathbb R \oplus\mathfrak n_{7,3,B}$ & $(e^{12},\,e^{23},\,e^{34})$ & $\begin{array}[t]{@{}l@{}} \Span( e^1+ie^2,\, e^3+i(e^4-e^1),\\ \phantom{\Span(}e^5+ie^6,\, e^7+ie^8) \end{array}$ \\
\hline
$\mathbb R \oplus\mathfrak n_{7,3,B_1}$ & $(e^{12}-e^{34},\,e^{13}+e^{24},\,e^{14})$ & $\begin{array}[t]{@{}l@{}} \Span( e^1+ie^4,\, e^2+ie^3,\\ \phantom{\Span(}e^5+ie^6,\, e^7+ie^8) \end{array}$ \\
\hline
$\mathbb R \oplus\mathfrak n_{7,3,C}$ & $(e^{12}+e^{34},\,e^{23},\,e^{24})$ & $\begin{array}[t]{@{}l@{}} \Span( e^1+ie^2,\, e^3+ie^4,\\ \phantom{\Span(}e^5+ie^6,\, e^7+ie^8) \end{array}$ \\
\hline
$\mathbb R \oplus\mathfrak n_{7,3,D}$ & $(e^{12}+e^{34},\,e^{13},\,e^{24})$ & $\begin{array}[t]{@{}l@{}} \Span( e^1+ie^3,\, e^2+ie^4,\\ \phantom{\Span(}e^5+ie^6,\, e^7+ie^8) \end{array}$ \\
\hline
$\mathbb R \oplus\mathfrak n_{7,3,D_1}$ & $(e^{12}-e^{34},\,e^{13}+e^{24},\,e^{14}-e^{23})$ & $\begin{array}[t]{@{}l@{}} \Span( e^1+ie^2,\, e^3+ie^4,\\ \phantom{\Span(}e^5+ie^6,\, e^7+ie^8) \end{array}$ \\
\hline
$\mathfrak n_{5,2}\oplus\mathfrak h_3$ & $(e^{12},\,e^{13},\,e^{45})$ & $\begin{array}[t]{@{}l@{}} \Span( e^1+ie^8,\, e^2+ie^3,\\ \phantom{\Span(}e^4+ie^5,\, e^6+ie^7) \end{array}$ \\
\hline
$\mathfrak n_{8,3}^{1}$ & $(e^{12}+e^{45},\,e^{25}+e^{34},\,e^{15}+e^{23})$ & $\begin{array}[t]{@{}l@{}} \Span( e^1-i(e^2+e^5),\, e^3+ie^5,\\ \phantom{\Span(}e^4+ie^8,\, e^6+ie^7) \end{array}$ \\
\hline
$\mathfrak n_{8,3}^{1\text{-bis}}$ & $(e^{12}+2e^{23}+2e^{45},\,e^{25}-2e^{34},\,-e^{15}+2e^{24}-2e^{35})$ & $\begin{array}[t]{@{}l@{}} \Span( e^1+ie^7,\, e^2+ie^5,\\ \phantom{\Span(}e^3-ie^4,\, e^6-ie^8) \end{array}$ \\
\hline
$\mathfrak n_{8,3}^{2}$ & $(e^{14}+e^{35},\,e^{15}+e^{24},\,e^{23})$ & $\begin{array}[t]{@{}l@{}} \Span( e^1+ie^4,\, e^2+ie^5,\\ \phantom{\Span(}e^3+ie^7,\, e^6+ie^8) \end{array}$ \\
\hline
$\mathfrak n_{8,3}^{2\text{-bis}}$ & $(e^{14}+e^{24}-e^{35},\,-e^{15}+e^{25}+e^{34},\,-e^{23})$ & $\begin{array}[t]{@{}l@{}} \Span( e^1+ie^8,\, e^2+ie^3,\\ \phantom{\Span(}e^4-ie^5,\, e^6+ie^7) \end{array}$ \\
\hline
$\mathfrak n_{8,3}^{3}$ & $(e^{15}+e^{34},\,e^{23},\,e^{14}+e^{24})$ & $\begin{array}[t]{@{}l@{}} \Span( e^1+ie^7,\, e^2-ie^3,\\ \phantom{\Span(}e^4-ie^5,\, e^6+ie^8) \end{array}$ \\
\hline
$\mathfrak n_{8,3}^{3\text{-bis}}$ & $(e^{15}+e^{34},\,e^{14}-e^{23},\,e^{13}+e^{24})$ & $\begin{array}[t]{@{}l@{}} \Span( e^1+ie^5,\, e^2+ie^6,\\ \phantom{\Span(}e^3+ie^4,\, e^7-ie^8) \end{array}$ \\
\hline
$\mathfrak n_{8,3}^{4}$ & $(e^{25}+e^{34},\,e^{15}+e^{23},\,e^{12})$ & $\begin{array}[t]{@{}l@{}} \Span( e^1-ie^2,\, e^3+i(e^1-e^4),\\ \phantom{\Span(}e^5+ie^8,\, e^6+ie^7) \end{array}$ \\
\hline
$\mathfrak n_{8,3}^{5\text{-bis}}$ & $(e^{13},\,e^{15}-e^{24},\,e^{14}+e^{25})$ & $\begin{array}[t]{@{}l@{}} \Span( e^1+ie^3,\, e^2+ie^6,\\ \phantom{\Span(}e^4-ie^5,\, e^7+ie^8) \end{array}$ \\
\hline
$\mathfrak n_{8,3}^{8}$ & $(e^{25}+e^{34},\,e^{12},\,e^{15})$ & $\begin{array}[t]{@{}l@{}} \Span( e^1-ie^2,\, e^3+ie^4,\\ \phantom{\Span(}e^5+ie^7,\, e^6+ie^8) \end{array}$ \\
\hline
$\mathfrak n_{8,3}^{9}$ & $(e^{13}+e^{25},\,e^{15}+e^{24},\,e^{12})$ & $\begin{array}[t]{@{}l@{}} \Span( e^1+ie^2,\, e^3+i(e^5+e^8),\\ \phantom{\Span(}e^4+i(e^8-e^5),\, e^6+ie^7)\end{array}$ \\
\hline
$\mathfrak n_{8,3}^{10}$ & $(e^{14}+e^{35},\,e^{23},\,e^{12})$ & $\begin{array}[t]{@{}l@{}} \Span( e^1+ie^3,\, e^2+ie^6,\\ \phantom{\Span(}e^4+ie^5,\, e^7+ie^8) \end{array}$ \\
\hline
$\mathfrak n_{8,3}^{12}$ & $(e^{15}+e^{23},\,e^{14},\,e^{12})$ & $\begin{array}[t]{@{}l@{}} \Span( e^1+ie^2,\, e^4+ie^8,\\ \phantom{\Span(}e^5+i(e^3+e^4),\, e^6+ie^7) \end{array}$ \\
\hline
\end{longtable}
}

{\small
\renewcommand{\arraystretch}{1.3}
\begin{longtable}{|l|l|l|}
\hline
$\n$ & $(de^5,de^6,de^7,de^8)$ & $\Lambda^{1,0}\n^*$ \\
\hline
\endfirsthead
\hline
$\n$ & $(de^5,de^6,de^7,de^8)$ & $\Lambda^{1,0}\n^*$ \\
\hline
\endhead
\multicolumn{3}{r}{\footnotesize\emph{Continued on next page}} \\
\endfoot
\caption{Complex structures on $8$-dimensional $2$-step nilpotent Lie algebras with $\dim(\n')=4$}
\label{tab:complex-4}
\endlastfoot
$\mathfrak n_{8,4}^{1}$ & $(e^{12},\,e^{34},\,e^{13},\,e^{24})$ & $\Span( e^1+ie^4,\, e^2-ie^3,\, e^5-ie^7,\, e^6+ie^8)$ \\
\hline
$\mathfrak n_{8,4}^{1\text{-bis}}$ & $(e^{12},\,e^{34},\,e^{13}-e^{24},\,-e^{14}-e^{23})$ & $\Span( e^1+ie^2,\, e^3-ie^4,\, e^5+ie^6,\, e^7+ie^8)$ \\
\hline
$\mathfrak n_{8,4}^{2}$ & $(e^{13},\,e^{24},\,e^{14}+e^{23},\,e^{12})$ & $\Span( e^1+ie^3,\, e^2-ie^4,\, e^5+ie^6,\, e^7+2ie^8)$ \\
\hline
$\mathfrak n_{8,4}^{3}$ & $(e^{14},\,e^{23},\,e^{12},\,e^{13})$ & $\Span( e^1+ie^4,\, e^2+ie^3,\, e^5+ie^6,\, e^7+ie^8)$ \\
\hline
\end{longtable}
}

\medskip
\section{Obstructions to integrability in dimension $8$}
\label{sec:5}

We consider in this section the case of $2$-step nilpotent $8$-dimensional Lie algebras $\n$ with center $\z$ of dimension $3$, equal to the commutator ideal $\n'$.

The aim of this section is to derive obstructions to the existence of complex structures on $\n$ in terms of the properties of its pencil $P_\n$. We start with the following:

\begin{lem} Assume that $J$ is a complex structure on $\n$ and denote by $\z_0:=\z\cap J\z$.
    Then $\dim(\z_0)=2$.
\end{lem}
\begin{proof}
    Since $\z_0$ is $J$-invariant, it has even dimension. Moreover, $\dim(\z)=3$ so the only possibilities are $\dim(\z_0)=0$ or $\dim(\z_0)=2$. We need to exclude the first case. 

    Assume that $\z\cap J\z=0$. Using that $N_J=0$, for every $x,y\in\n$ we obtain 
    $$[Jx,Jy]-[x,y]=J([x,Jy]+[Jx,y])\in \z\cap J\z=0,$$
    so $[Jx,Jy]=[x,y]$. For every $x\in \z$, this implies $Jx\in\z$, which contradicts the assumption.    
\end{proof}

Let us fix from now on a Hermitian metric $g$ on $(\n,J)$.
Let $\{z_1,z_2,z_3\}$ be an orthonormal basis of $\z$ with $z_2=Jz_1\in\z_0$. We denote by $\v:=\z^\perp$ and by $\xi:=Jz_3\in \v$. Let $\v_0$ be the orthogonal of $\xi$ in $\v$, so that we have the following $J$-invariant orthogonal decomposition of $\n$:
$$\n=\v_0\oplus\Span(\xi,z_3)\oplus\z_0.$$

We denote by $z^i=g(z_i,\cdot)$ for $i=1,2,3$, and $\xi^\flat:=g(\xi,\cdot)$ the corresponding metric duals in $\n^*$. The differentials $dz^i$ vanish on $\z$ so can be seen as elements of $\Lambda^2\v^*=\Lambda^2\v_0^*\oplus(\xi^\flat\wedge \v_0^*).$ Using also the direct sum decomposition $\Lambda^2\v_0^*=L^+(\v_0)\oplus L^-(\v_0)$ introduced in \eqref{lambda}, we write the 2-forms $dz^i$ for $i=1,2,3$ as
\begin{equation}\label{dzi}
    dz^i=\xi^\flat\wedge\eta_i+\alpha_i+\beta_i,\qquad\hbox{with } \eta_i\in\v_0^*,\ \alpha_i\in L^+(\v_0),\ \beta_i\in L^-(\v_0).
\end{equation}
\begin{lem}\label{l42}
    The structure $J$ is integrable if and only if 
    \begin{equation}\label{int}
        \eta_3=0,\qquad\beta_3=0,\qquad \eta_2=-\eta_1\circ J,\qquad \beta_2=-\beta_1(J\cdot,\cdot).
    \end{equation}
\end{lem}
\begin{proof}
  Using \eqref{nij} and the fact that $N_J(x,y)\in\z+J\z=\Span(z_1,Jz_1,z_3,Jz_3)$ for every $x,y\in\n$,
    the vanishing of $N_J$ is equivalent to 
    $$0=g(N_J(x,y),z_1)=g(N_J(x,y),z_3),\qquad\forall x,y\in\n.$$
    Since $Jz_3=\xi$ is orthogonal to $\n'$, these two equations are equivalent to the system 
    \begin{equation}\label{sys}
        \begin{cases}
            dz^3(x,y)=dz^3(Jx,Jy),\qquad&\forall x,y\in\n,\\
            dz^1(x,y)-dz^1(Jx,Jy)=dz^2(Jx,y)+dz^2(x,Jy),\qquad&\forall x,y\in\n.
        \end{cases}
    \end{equation}
     Using \eqref{dzi}, the system \eqref{sys} for $x=\xi$ is equivalent to 
    \begin{equation}\label{sys1}
        \begin{cases}
            \eta_3=0,\\
            \eta_1=\eta_2\circ J,
        \end{cases}
    \end{equation}
    whereas for $x,y\in\v_0$, \eqref{sys} is equivalent to 
    \begin{equation}\label{sys2}
        \begin{cases}
            \beta_3=0,\\
            \beta_1=\beta_2(J\cdot,\cdot).
        \end{cases}
    \end{equation}
    Since $N_J(x,y)=0$ whenever $x$ or $y$ belong to $\z_0$, this concludes the proof.
\end{proof}

    Note that $\eta_1\neq0$, since otherwise $\xi$ would belong to the center of $\n$, contradicting the assumption $\dim(\z)=3$.

In view of Lemma \ref{l42}, and Remark \ref{rem2}, we make the following
\begin{defi}
    Let $\v$ be a $5$-dimensional real vector space. A $3$-dimensional subspace $P\subset\Lambda^2\v^*$ is called a pencil of \textit{complex type} if there exists an inner product $g$ on $\v$, a unit vector $\xi\in\v$ with dual $1$-form $\xi^\flat:=g(\xi,\cdot)$, a Hermitian structure $J$ on $\v_0:=\xi^\perp$, and elements $\eta\in\v_0^*\setminus\{0\},\ \alpha_1,\alpha_2,\alpha_3\in L^+(\v_0),\ \beta\in L^-(\v_0)$ such that the 2-forms $\omega_1,\ \omega_2,\ \omega_3$ defined by
    \begin{equation}\label{sys3}
    \begin{cases}
        \omega_1:= \xi^\flat\wedge\eta\ \hskip 0,2cm+&\hskip -0,2cm\alpha_1+\beta,\\
        \omega_2:=\xi^\flat\wedge J\eta\ +&\hskip -0,2cm\alpha_2+J\circ\beta,\\
       \omega_3:=&\hskip -0,2cm\alpha_3,\\
    \end{cases}
\end{equation}
are a basis of $P$.
\end{defi}

\begin{rem}\label{re}
    By replacing $\eta$ with $\cos \theta\,\eta+\sin \theta\, J\eta$, $\beta$ with $\cos \theta\,\beta+\sin \theta\, J\circ\beta$, $\alpha_1$ with $\cos \theta\,\alpha_1+\sin \theta \,\alpha_2$, and $\alpha_2$ with $-\sin \theta\,\alpha_1+\cos \theta\, \alpha_2$, the above basis $\{\omega_1,\omega_2,\omega_3\}$ of $P$ is replaced by
    $$\{\tilde\omega_1,\tilde\omega_2,\tilde\omega_3\}:=\{\cos \theta\,\omega_1+\sin \theta \,\omega_2, -\sin \theta\,\omega_1+\cos \theta \,\omega_2,\omega_3\}.$$
    In particular, after such a change of basis in $P$, every element $\gamma\in P\setminus\Span(\omega_3)$ can be written as $\rho(\omega_1+t\omega_3)$ for some $\rho\in\R\setminus\{0\}$. Indeed, if $\gamma=a\,\omega_1+b\,\omega_2+c\,\omega_3$ with $(a,b)\neq(0,0)$, then writing $a+ib=\rho e^{i\theta}$ gives $\gamma=\rho(\tilde\omega_1+t\tilde\omega_3)$, for $t:=\frac c\rho$.
\end{rem}

\begin{cor}\label{c}
    If $\n$ admits a complex structure, its pencil $P_\n\subset \Lambda^2\v^*$ is of complex type.
\end{cor}
\begin{proof}
    Follows directly from Lemma \ref{l42}, together with Equation \eqref{dzi}, by denoting $\eta:=\eta_1$, and $\beta:=\beta_1$.
\end{proof}

  Using Corollary \ref{c}, we will now obtain obstructions for certain $2$-step nilpotent $8$-dimen\-sional Lie algebras with $3$-dimensional center to admitting complex structures.

  \begin{lem}\label{l34}
      Let $P\subset\Lambda^2\v^*$ be a pencil of complex type. Then the set 
      $$S(P):=\{\gamma\in P\ |\ \gamma\wedge\gamma=0\},$$
    of decomposable $2$-forms in $P$, is equal to the union of at most two vector spaces.
  \end{lem}
  \begin{proof}
  Let $\{\omega_1,\omega_2,\omega_3\}$ be a basis of $P$ satisfying \eqref{sys3}.
      Consider any element $\gamma=\sum_{i=1}^3t_i\omega_i\in P$. By \eqref{sys3}, one can write 
      $$\gamma=\xi^\flat\wedge(t_1\eta+t_2J\eta)+\sum_{i=1}^3t_i\alpha_i+t_1\beta+t_2 J\circ\beta.$$
     Since $\Lambda^4\v^*=\Lambda^4\v_0^*\oplus (\Lambda^3\v_0^*\wedge \xi^\flat)$, we see that $\gamma$ belongs to $S(P)$ if and only if 
      \begin{equation}\label{sp}
          \begin{cases}
              (t_1\eta+t_2J\eta)\wedge(\sum_{i=1}^3t_i\alpha_i+t_1\beta+t_2 J\circ\beta)=0,\\
              (\sum_{i=1}^3t_i\alpha_i+t_1\beta+t_2 J\circ\beta)^2=0.
          \end{cases}
      \end{equation}
      Note that if $t_1=t_2=0$, the above system is equivalent to $t_3^2\alpha_3\wedge\alpha_3=0$, so 
       $$S(P)\cap \Span(\omega_3)=\begin{cases}
          0&\hbox{ if }\alpha_3\wedge\alpha_3\neq0,\\
          \Span(\omega_3)&\hbox{ if }\alpha_3\wedge\alpha_3=0.
      \end{cases}$$
     On the other hand, if 
      $(t_1,t_2)\neq(0,0)$, then $t_1\eta+t_2J\eta\neq 0$, so the first equation of \eqref{sp} implies the second one. Thus, using also \eqref{eb},  
      $S(P)\setminus \Span(\omega_3)$ is the set of all $\gamma=\sum_{i=1}^3t_i\omega_i$ such that $(t_1,t_2)\neq(0,0)$ and 
      \begin{equation}\label{e9}
          (t_1\eta+t_2J\eta)\wedge\sum_{i=1}^3t_i\alpha_i+(t_1^2+t_2^2)\eta\wedge\beta=0.
      \end{equation}
      Assume that \eqref{e9} holds. Applying $J$ as an automorphism of $\Lambda\v_0^*$ (see Section \ref{sec23}) yields 
       \begin{equation}\label{e10}
          (t_1J\eta-t_2\eta)\wedge\sum_{i=1}^3t_i\alpha_i-(t_1^2+t_2^2)J\eta\wedge\beta=0.
      \end{equation}
      Multiplying \eqref{e9} by $t_1$, \eqref{e10} by $t_2$ and subtracting yields 
      \begin{equation}\label{e11}
          (t_1^2+t_2^2)\eta\wedge\sum_{i=1}^3t_i\alpha_i+(t_1^2+t_2^2)(t_1\eta+t_2J\eta)\wedge\beta=0,
      \end{equation}
      and since $(t_1,t_2)\neq(0,0)$, this is equivalent to 
       \begin{equation}\label{e12}
        \eta\wedge\sum_{i=1}^3t_i\alpha_i+(t_1\eta+t_2J\eta)\wedge\beta=0.
      \end{equation}
      Conversely, applying $J$ to \eqref{e12} and taking a suitable linear combination gives \eqref{e9}.
      Let us denote by $Q\subset P$ the set of solutions of \eqref{e12}. Clearly $Q$ is a vector subspace of $P$. The above argument shows that $S(P)\setminus \Span(\omega_3)=Q\setminus \Span(\omega_3)$, whence
       \begin{equation}\label{spq}
          S(P)=\left(S(P)\setminus \Span(\omega_3)\right)\cup\left(S(P)\cap\Span(\omega_3)\right)=\left(Q\setminus \Span(\omega_3)\right)\cup\left(S(P)\cap\Span(\omega_3)\right). 
          \end{equation}
      Now, obviously 
      $$Q\cap \Span(\omega_3)=\begin{cases}
          0&\hbox{ if }\eta\wedge\alpha_3\neq0\\
          \Span(\omega_3)&\hbox{ if }\eta\wedge\alpha_3=0.
      \end{cases}$$
      We have three cases:
      \begin{itemize}
          \item if $\eta\wedge\alpha_3= 0$, then also $\alpha_3\wedge\alpha_3=0$, so 
          $$S(P)\cap\Span(\omega_3) =\Span(\omega_3)=Q\cap\Span(\omega_3)
         ,$$
      whence by \eqref{spq}, $S(P)=\left(Q\setminus \Span(\omega_3)\right)\cup (Q\cap\Span(\omega_3))=Q$.
      \item if $\eta\wedge\alpha_3\neq0$, but $\alpha_3\wedge\alpha_3=0$, then
          $$S(P)\cap\Span(\omega_3)=\Span(\omega_3),\qquad Q\cap\Span(\omega_3)=0,$$
      whence by \eqref{spq} we get $S(P)=(Q\setminus\{0\})\cup\Span(\omega_3)=Q\cup\Span(\omega_3).$
      \item if $\alpha_3\wedge\alpha_3\neq0$, then $\eta\wedge\alpha_3\neq0$, so 
       $$S(P)\cap\Span(\omega_3)=Q\cap\Span(\omega_3)
          =0,$$
      whence by \eqref{spq} again, $S(P)=Q$.
      \end{itemize}
     This proves the lemma.
  \end{proof}

  \begin{lem}\label{l35}
       Let $P\subset\Lambda^2\v^*$ be a pencil of complex type. If there exist $\gamma,\gamma'\in P\setminus\Span(\omega_3)$ with $\omega_3\wedge\gamma=0$ and $\gamma'\wedge\gamma'=0$, then the set $S(P)$, of decomposable $2$-forms in $P$, contains a $2$-plane.
  \end{lem}
  \begin{proof}
     Let $\{\omega_1,\omega_2,\omega_3\}$ be a basis of $P$ satisfying \eqref{sys3}. Denote 
      $$\gamma=\sum_{i=1}^3t_i\omega_i,\qquad \gamma'=\sum_{i=1}^3t'_i\omega_i$$
      with $(t_1,t_2)\neq(0,0)\neq(t'_1,t'_2)$. Projecting onto the second summand of $\Lambda^4\v^*=\Lambda^4\v_0^*\oplus (\Lambda^3\v_0^*\wedge \xi^\flat)$, the first condition implies $\alpha_3\wedge(t_1\eta+t_2J\eta)=0$. Applying $J$ we obtain $\alpha_3\wedge(t_1J\eta-t_2\eta)=0$, so $\alpha_3\wedge\eta=\alpha_3\wedge J\eta=0$. Thus $\omega_3=\alpha_3=c\eta\wedge J\eta$ for some $c\neq 0$.

      By Remark \ref{re}, one can assume that $$\gamma'=\omega_1+t\omega_3=\xi^\flat\wedge\eta+(\alpha_1+t\alpha_3+\beta).$$ Then $\gamma'\wedge\gamma'=0$ implies $\eta\wedge(\alpha_1+t\alpha_3+\beta)=0$, so $\eta\wedge\gamma'=0$, whence $\Span(\omega_3,\gamma')\subset S(P)$.
  \end{proof}

\begin{lem}\label{l36}
       Let $P\subset\Lambda^2\v^*$ be a pencil of complex type. If there exist $\gamma,\gamma'\in P\setminus\Span(\omega_3)$ with $\omega_3\wedge\gamma=0$ and $\gamma\wedge\gamma=\omega_3\wedge\gamma'$, then $S(P)$ contains a $2$-plane.
  \end{lem}
  \begin{proof}
     Let $\{\omega_1,\omega_2,\omega_3\}$ be a basis of $P$ satisfying \eqref{sys3}. Denote 
      $$\gamma=\sum_{i=1}^3t_i\omega_i,\qquad \gamma'=\sum_{i=1}^3t'_i\omega_i$$
      with $(t_1,t_2)\neq(0,0)\neq(t'_1,t'_2)$. The argument in the previous lemma shows that $\omega_3=\alpha_3=c\eta\wedge J\eta$ for some $c\neq 0$.
      Using again Remark \ref{re}, one can assume this time that $$\gamma=\omega_1+t\omega_3=\xi^\flat\wedge\eta+(\alpha_1+t\alpha_3+\beta).$$ Then projecting $\gamma\wedge\gamma=\omega_3\wedge\gamma'$ onto $\Lambda^3\v_0^*\wedge \xi^\flat$ implies
      $$2\eta\wedge(\alpha_1+t\alpha_3+\beta)=\alpha_3\wedge(t'_1\eta+t'_2J\eta)=0.$$
     We thus obtain $\eta\wedge\gamma=0$, whence $\Span(\omega_3,\gamma)\subset S(P)$.
  \end{proof}
  Finally, we prove the following
  \begin{lem}\label{l37}
       Let $P\subset\Lambda^2\v^*$ be a pencil of complex type. If there exist $\gamma\in P\setminus\Span(\omega_3)$ and $\gamma'\in P\setminus\Span(\gamma)$ with $\gamma\wedge\gamma=\gamma\wedge\gamma'=0$, then $S(P)$ contains a $2$-plane.
  \end{lem}
  \begin{proof}
    Let $\{\omega_1,\omega_2,\omega_3\}$ be a basis of $P$ satisfying \eqref{sys3}.  By Remark \ref{re} again, one can assume that $$\gamma=\omega_1+t\omega_3=\xi^\flat\wedge\eta+(\alpha_1+t\alpha_3+\beta).$$
      Since $\gamma\wedge\gamma=0$, we obtain as before $\eta\wedge(\alpha_1+t\alpha_3+\beta)=0$. Applying $J$, this also gives
      \begin{equation}\label{13}
          J\eta\wedge(\alpha_1+t\alpha_3-\beta)=0.
      \end{equation}
      Consider now
      $$\gamma'=\sum_{i=1}^3t'_i\omega_i.$$
      Denoting $\gamma'':=\gamma'-t'_1\gamma$, one thus has $\gamma''=s_2\omega_2+s_3\omega_3$, and
     $$ 
0=\gamma\wedge\gamma''=(\xi^\flat\wedge\eta+(\alpha_1+t\alpha_3+\beta))\wedge(s_2\xi^\flat\wedge J\eta+s_2\alpha_2+s_2J\circ\beta+s_3\alpha_3).
    $$ 
Projecting again onto $\Lambda^3\v_0^*\wedge \xi^\flat$ and using \eqref{13} and then \eqref{eb} yields
\begin{eqnarray*}
0&=&\eta\wedge(s_2\alpha_2+s_2J\circ\beta+s_3\alpha_3)+s_2J\eta\wedge(\alpha_1+t\alpha_3+\beta)\\
&=&\eta\wedge(s_2\alpha_2+s_2J\circ\beta+s_3\alpha_3)+2s_2J\eta\wedge\beta\\
&=&\eta\wedge(s_2\alpha_2-s_2J\circ\beta+s_3\alpha_3).
\end{eqnarray*}
Applying $J$ to this equation we obtain
$$0=J\eta\wedge(s_2\alpha_2+s_2J\circ\beta+
s_3\alpha_3)=J\eta\wedge\gamma''.$$
 Therefore $\gamma''\wedge\gamma''=0$. Moreover, $\gamma\wedge\gamma''=\gamma\wedge\gamma'=0$, and $\gamma''$ is not proportional to $\gamma$ (since $\gamma'\in P\setminus\Span(\gamma)$), whence $\Span(\gamma,\gamma'')$ is a 2-plane contained in $ S(P)$.
  \end{proof}

  \medskip
\section{$2$-step nilpotent Lie algebras with no complex structures} \label{sec:6}

Using the obstructions found in the previous section, we will now finish the proof of the ``only if" part of Theorem \ref{thmB}, by showing that the Lie algebras $\n^5_{8,3}$, $\n^6_{8,3}$, and $\n^7_{8,3}$ (see the Appendix for their definition) carry no complex structures. These are the first examples, to the best of our knowledge, of $2$-step nilpotent Lie algebras with no complex structure. They are all $8$-dimensional, with center and commutator ideal of dimension $3$.

\subsection{The algebra $\n^5_{8,3}$}\label{n538} A pencil of $\n^5_{8,3}$ is $P_{\n^5_{8,3}}:=\Span(e^{13},e^{24},e^{15}+e^{25})$. We compute 
$$\frac12\left(a\,e^{13}+b\,e^{24}+c\,(e^{15}+e^{25})\right)^{\wedge2}=-ab\,e^{1234}-ac\,e^{1235}+bc\,e^{1245},$$
so $S(P_{\n^5_{8,3}})=\Span(e^{13})\cup\Span(e^{24})
\cup\Span(e^{15}+e^{25})$. By Lemma \ref{l34}, $P_{\n^5_{8,3}}$ is not of complex type, so by Corollary \ref{c}, $\n^5_{8,3}$ does not admit complex structures.

\subsection{The algebra $\n^6_{8,3}$}\label{n683} A pencil of $\n^6_{8,3}$ is $P_{\n^6_{8,3}}:=\Span(e^{24},e^{15}+e^{23},e^{13})$. We compute 
$$\frac12\left(a\,e^{24}+b\,(e^{15}+e^{23})+c\,e^{13})\right)^{\wedge2}=ab\,e^{1245}-ac\,e^{1234}+b^2\,e^{1235},$$
so $S(P_{\n^6_{8,3}})=\Span(e^{13})\cup\Span(e^{24})$. Assume that $P_{\n^6_{8,3}}$ is of complex type, so satisfies the system \eqref{sys3}. 

We distinguish two cases. If $e^{13}\in P_{\n^6_{8,3}}\setminus\Span(\omega_3)$, this contradicts Lemma \ref{l37} applied to $\gamma:=e^{13}$ and $\gamma':=e^{15}+e^{23}$. On the other hand, if $e^{13}\in\Span(\omega_3)$, this contradicts Lemma \ref{l35} applied to $\gamma:=e^{15}+e^{23}$ and $\gamma':=e^{24}$.  

Thus $P_{\n^6_{8,3}}$ is not of complex type, so by Corollary \ref{c}, $\n^6_{8,3}$ does not admit complex structures.

\subsection{The algebra $\n^7_{8,3}$} A pencil of $\n^7_{8,3}$ is $P_{\n^7_{8,3}}:=
\Span(e^{15}+e^{24},e^{14}+e^{23},e^{13})$. We compute as before 
$$\frac12\left(a\,(e^{15}+e^{24})+b\,(e^{14}+e^{23})+c\,e^{13})\right)^{\wedge2}=a^2\,e^{1245}+ab\,e^{1235}+(b^2-ac)\,e^{1234},$$
so $S(P_{\n^7_{8,3}})=\Span(e^{13})$. Assume that $P_{\n^7_{8,3}}$ is of complex type, so satisfies the system \eqref{sys3}.

We distinguish again two cases. If $e^{13}\in P_{\n^7_{8,3}}\setminus\Span(\omega_3)$, this contradicts Lemma \ref{l37} applied to $\gamma:=e^{13}$ and $\gamma':=e^{14}+e^{23}$. On the other hand, if $e^{13}=c\,\omega_3$ for some non-zero $c$, this contradicts Lemma \ref{l36} applied to $\gamma:=e^{14}+e^{23}$ and $\gamma':=-2c\,(e^{15}+e^{24})$.  

Thus $P_{\n^7_{8,3}}$ is not of complex type, so by Corollary \ref{c}, $\n^7_{8,3}$ does not admit complex structures.

\medskip
\subsection{Examples in higher dimensions} 

By applying repeatedly Theorems~\ref{thm:g+r2} and~\ref{thm:g+r2'} and Proposition~\ref{prop34} to $\n_{8,3}^{i}$, for $i=5,6,7$, we obtain, in every even dimension $\geq10$, both decomposable and indecomposable examples of $2$-step nilpotent Lie algebras which do not carry any complex structure, thus completing the proof of Theorem \ref{thmB}.

\begin{thm}
Let $i\in\{5,6,7\}$, $m\geq1$ and $\xi\in\z(\n_{8,3}^{i})$. Then the $(2m+8)$-dimensional $2$-step nilpotent Lie algebra
\[
\n^{i,\varphi_m}_{2m+8,3}:=\R^{2m}\oplus_{\varphi_m}\n_{8,3}^{i},\qquad
\varphi_m:=\big(e^{12}+e^{34}+\cdots+e^{2m-1,2m}\big)\otimes\xi ,
\]
admits no complex structure.
For $\xi=0$,  $\n^{i,\varphi_m}_{2m+8,3}=\R^{2m}\oplus\n_{8,3}^{i}$ is decomposable, while for $\xi\neq0$, $\n^{i,\varphi_m}_{2m+8,3}$ is indecomposable.
\end{thm}
\begin{proof}
      The case $\xi=0$ follows directly by applying $m$ times Theorem \ref{thm:g+r2}. For $\xi\neq 0$, using the isomorphism
$$\n^{i,\varphi_m}_{2m+8,3}\simeq\R^2\oplus_{\varphi_1}\n^{i,\varphi_{m-1}}_{2m+6,3},$$ the non-existence of complex structures follows by applying  $m$ times Theorem \ref{thm:g+r2'}, and the indecomposability follows by induction on $m$ from Proposition \ref{prop34}.
\end{proof}

 \subsection{Complex structures on extensions of $\R^2$ by Lie algebras without complex structures}\label{s65}
     Consider the $10$-dimensional 2-step nilpotent Lie algebra $\n:=\Span(e_1,\ldots,e_7,z_1,z_2,z_3)$, whose differentials of the dual basis are given by:
$$dz^1=e^{24},\qquad dz^2=e^{15}+e^{23}+e^{67}, \qquad dz^3=e^{13}+e^{57}+e^{46},\qquad de^i=0,\ \forall 1\leq i\leq 7.$$
Define the almost complex structure $J$ on $\n$ whose $(1,0)$ space $\Lambda^{1,0}\n^*$ is spanned by 
$$\sigma_1:=e^1+ie^5,\qquad \sigma_2:=e^2+ie^6,\qquad \sigma_3:=e^3+ie^7,\qquad \sigma_4:=e^4-iz^2,\qquad \sigma_5:=z^1-iz^3.$$
We have 
$$d\sigma_1=d\sigma_2=d\sigma_3=0,\qquad d\sigma_4=-i(e^{15}+e^{23}+e^{67})=\frac12\sigma_1\wedge\bar\sigma_1-i \mathrm{Re}(\sigma_2\wedge\bar\sigma_3)\in\Lambda^{1,1}\n^*,\text{ and}$$
$$d\sigma_5=e^{24}-i(e^{13}+e^{57}+e^{46})=-i\mathrm{Re}(\sigma_1\wedge\bar\sigma_3)+\frac{1}2\sigma_2\wedge(\sigma_4+\bar\sigma_4)\in\Lambda^{1,1}\n^*\oplus\Lambda^{2,0}\n^*,$$
thus showing that $J$ is integrable. On the other hand, the ideal $\n_0:=\Span(e_1,\ldots,e_5,z_1,z_2,z_3)$ of $\n$ is determined by the differentials 
$$dz^1=e^{24},\qquad dz^2=e^{15}+e^{23}, \qquad dz^3=e^{13},\qquad de^i=0,\ \forall 1\leq i\leq 5.$$
This shows that $\n_0$ coincides with  $\n^6_{8,3}$, which admits no complex structures by \S\ref{n683}. Note that $\n$ is an extension of $\R^2$ by $\n_0$ since there is an obvious exact sequence
$$0\to\n_0\to\n\to\R^2\to 0.$$

\medskip

\section{Complex structures on $2$-step nilpotent Lie algebras with $2$-dimensional center}\label{sec:7}

Every $2$-step nilpotent Lie algebra $\n$ with $2$-dimensional center $\z$ also has $2$-dimensional commutator ideal $\n'$, except for $\R\oplus\h_{2n-1}$, which obviously admits abelian complex structures. Indeed, any complex structure on   $\R\oplus\h_{2n-1}$ is abelian \cite[Proposition 3.6]{Rol}. We will thus concentrate on those for which $\n'=\z$.

The set of isomorphism classes of $2$-step nilpotent Lie algebras $\n$ satisfying $\dim \z =2$, $\n '=\z$, can be  parametrized using pencil invariants (see \cite{Gau,LT}). In this section we  recall from \cite{Lau-Osc} the description of such Lie algebras and show that they always admit complex structures, which  in many cases are abelian (Theorem \ref{thm:dimz=2} below).

 By \cite[Theorems~3.6 and~4.4]{LT}
(see \cite[Section 6]{Gau} for the complex case), the isomorphism classes of such
Lie algebras over $\R$ are parametrized by sets of the form
\begin{equation}\label{eq:parameters}
  \mathcal S=\{(\alpha_1,k_1), \ldots ,(\alpha_r,k_r), \epsilon_1, \ldots, \epsilon_s\},
  \qquad \alpha_\ell\in \C\cup \{\infty\},\; k_\ell,\,\epsilon_j\in \N .
\end{equation}
We recall this correspondence below, following the lines of \cite[Section 3]{Lau-Osc}. Consider the vector spaces 
  \begin{equation} \label{eq:nS}     
  \v:=\v_1 \oplus \cdots \oplus \v_r\oplus\fu_1 \oplus \cdots \oplus \fu_s,\qquad\n_{\mathcal S}:=\z\oplus\v ,
  \end{equation}
  with $\dim \z=2$, $\dim \v_\ell =\begin{cases} 2k_\ell ,& \text{if  } \alpha_\ell\in \R\cup \{\infty\},\\  4k_\ell , & \text{if } \alpha_\ell\in \C\setminus\R , \end{cases}$ and $\dim \fu_j=2\epsilon_j+1 $. 
  Fix an inner product on $\n_\mathcal S$ such that the  decompositions \eqref{eq:nS} are orthogonal. 
    Let  
 $\{z_1, z_2\}$ be an orthonormal basis of $\z$ and fix 
  orthonormal bases of $\v_\ell$ and $\fu_j$. The Lie bracket on $\n_\mathcal S$  is determined by the skew-symmetric endomorphisms $j(z_1), \, j(z_2)\in \mathfrak{so}(\v)$ defined in Section \ref{sec:prelim}.  Their
   restrictions  to the subspaces $\v_\ell$ with respect to the chosen bases are given below.

   When $\alpha_\ell =\infty$, then:
   \begin{align}\label{eq:jz1}
   j(z_1)|_{\v_\ell} 
   &=\begin{bmatrix} 0& P_\ell \\
   -P_\ell & 0
   \end{bmatrix}, \quad  \text{ with }
P_\ell=\begin{bmatrix} & & 1\\ & \reflectbox{$\ddots$} & \\ 1 & & \end{bmatrix}, \\   \label{eq:jz2}
   j(z_2)\mid_{\v_\ell} &= \begin{bmatrix} 0& C_\ell \\
   -C_\ell & 0
   \end{bmatrix},\quad  \text{ with } C_\ell=\begin{bmatrix}
 & & & 0\\
 & & 0 & 1\\
 & \reflectbox{$\ddots$} & \reflectbox{$\ddots$} & \\
0 & 1 & &
\end{bmatrix}. \end{align}
 When $\alpha_\ell =a_\ell\in \R $, then 
 \begin{align}\label{eq:jz3}
   j(z_1)\mid_{\v_\ell} &= \begin{bmatrix} 0& C_\ell \\
   -C_\ell & 0
   \end{bmatrix},\quad  \text{ with } C_\ell=\begin{bmatrix}
 & & & a_\ell\\
 & & a_\ell & 1\\
 & \reflectbox{$\ddots$} & \reflectbox{$\ddots$} & \\
a_\ell & 1 & &
\end{bmatrix},\\   \label{eq:jz4}
   j(z_2)|_{\v_\ell} 
   &=\begin{bmatrix} 0& P_\ell \\
   -P_\ell & 0
   \end{bmatrix}, \quad  \text{ with }
P_\ell=\begin{bmatrix} & & 1\\ & \reflectbox{$\ddots$} & \\ 1 & & \end{bmatrix}. \end{align}
 For $\alpha_\ell\in \R\cup\{ \infty\} $, the matrices $C_\ell$ and $P_\ell$ have size $k_\ell \times k_\ell$.
 
 When $\alpha_\ell=a_\ell+ib_\ell\in \C\setminus \R$, then 
  \begin{align}\label{eq:jz5}
   j(z_1)\mid_{\v_\ell} &= \begin{bmatrix} 0& C_\ell \\
   -C_\ell & 0
   \end{bmatrix},\quad  \text{ with } C_\ell=\begin{bmatrix}
       &&&&&b_\ell&a_\ell \\
           &&&&&a_\ell&-b_\ell \\
               &&&b_\ell&a_\ell&0& 1\\
                   &&&a_\ell&-b_\ell&1&0 \\
                       &&\reflectbox{$\ddots$} &&\reflectbox{$\ddots$}&& \\
                           b_\ell&a_\ell&0& 1 &&&\\
                               a_\ell&-b_\ell &1&0&&&
   \end{bmatrix} , \\ \label{eq:jz6}
   j(z_2)|_{\v_\ell} 
   &=\begin{bmatrix} 0& P_\ell \\
   -P_\ell & 0
   \end{bmatrix}, \quad  \text{ with }
P_\ell=\begin{bmatrix} & & 1\\ & \reflectbox{$\ddots$} & \\ 1 & & \end{bmatrix}  . \end{align}
 In this case, the matrices $C_\ell$ and $P_\ell$ have size $2k_\ell \times 2k_\ell$.

Furthermore, the restrictions of $j(z_1), \, j(z_2)\in \mathfrak{so}(\v)$ to $\fu_j$ for $1\leq j\leq s,$ are given, with respect to the chosen bases, by:
\begin{equation}  \label{eq:jz-odd-v}
   j(z_1)\mid_{\fu_j} 
   =\begin{bmatrix} 0& Q_j \\
   -Q_j^T & 0
   \end{bmatrix}, \qquad  \quad
   j(z_2)\mid_{\fu_j} = \begin{bmatrix} 0& R_j \\
   -R_j^T & 0
   \end{bmatrix}, 
   \end{equation}
where $Q_j$ and $R_j$ are the following $(\epsilon_j+1) \times \epsilon _j$ matrices:
\[
Q_j=\begin{bmatrix}
0 & \cdots & 0 \\
1 &        &   \\
  & \ddots &   \\
  &        & 1
\end{bmatrix},
\qquad \quad
R_j=\begin{bmatrix}
1 &        &   \\
  & \ddots &   \\
  &        & 1 \\
0 & \cdots & 0
\end{bmatrix}.
\]

Although this will not be needed in what follows, we point out that,  
as observed in  \cite[Section 3]{Lau-Osc}, two Lie algebras $\n_\mathcal S$ and $\n_{\mathcal S'}$
are isomorphic if and only if $r=r', \, s=s', \,
\{\epsilon_1, \ldots, \epsilon_s\}=\{\epsilon_1 ', \ldots, \epsilon_s'\}$ counting multiplicities, and there exists a real Möbius transformation $T:\C \cup \{\infty\} \to \C \cup \{\infty\}$,
\[Tz:= \frac{az+b}{cz+d} , \qquad\quad \begin{bmatrix}
   a&b \\ c&d 
\end{bmatrix}\in \mathrm{GL}(2, \R),
\]
such that 
\[ \mathcal \{(\alpha_1',k_1'), \ldots ,(\alpha_r',k_r')\}= \{(T\alpha_1,k_1), \ldots ,(T\alpha_r,k_r)\}, \quad \text{counting multiplicities}.\]
 
 Using the above description, 
 we show next that any even-dimensional $2$-step nilpotent Lie algebra with $2$-dimensional center admits a complex structure. Before stating this result, we recall from \cite[Definition 1.4]{Eb} that a  $2$-step nilpotent Lie algebra $\n =\z\oplus \v$ is called non-singular if for each non-zero $v \in \v $, the map $\v \to \z, \, u\mapsto [v, u]$ is surjective.  
\begin{thm}    
\label{thm:dimz=2} Any even-dimensional  $2$-step nilpotent Lie algebra  with $2$-dimensional center admits a complex structure.  

Moreover, the Lie algebra $\n_\mathcal S$ with $\mathcal S$ given by 
\begin{equation}\label{eq:S-abelian}    
\mathcal S=\{(\alpha_1,k_1), \ldots ,(\alpha_r,k_r)  \}, \quad \alpha_\ell\in \C \cup \{\infty \}, \; k_\ell \in \N ,\end{equation}  admits an abelian complex structure.  In particular, this holds for any   non-singular $2$-step nilpotent Lie algebra with $2$-dimensional center.
\end{thm}
\begin{proof}
As mentioned at the beginning of this section, in case $\dim\n'=1$ we have that $\n$ is isomorphic to $\R\oplus \h_{2n-1}$ and it follows from \cite[Proposition 3.6]{Rol} that $\n $ has a complex structure, which is necessarily abelian.

If $\dim\n'=2$ and $\n$ is even dimensional, then $\n$ is isomorphic to 
$\n_\mathcal S$ for some $\mathcal S$ as in \eqref{eq:parameters} with even $s=2m$. 
We show next that $\n_\mathcal S$ admits a complex structure. 
We have that $\n_\mathcal S=\z\oplus \v$, with $\v=\v_1 \oplus \cdots \oplus \v_r \oplus\fu_1\oplus \cdots \oplus \fu_s$. The endomorphisms $j(z_1)|_{\v_\ell}, \, j(z_2)|_{\v_\ell}$  are given by \eqref{eq:jz1} and \eqref{eq:jz2} if $\alpha_\ell=\infty$, \eqref{eq:jz3} and \eqref{eq:jz4} if $\alpha_\ell\in \R $, or   
\eqref{eq:jz5} and \eqref{eq:jz6} if $\alpha_\ell\in \C\setminus\R $. On the other hand, $j(z_1)|_{\fu_j}, \, j(z_2)|_{\fu_j}$  are given by \eqref{eq:jz-odd-v}.   
Set  $\fu:= {\fu_1}\oplus \cdots \oplus {\fu_s} $, we wish to define an endomorphism $J_\fu $ of $\fu$ such that   $\fu_{2i-1}\oplus\fu_{2i}$ is $J_\fu$-invariant for $1\leq i \leq m$. It suffices to define $J_\fu$ on $\fu_1\oplus \fu_2$, with analogous definitions on the remaining spaces $\fu_{2i-1}\oplus\fu_{2i}$, $i\geq 2$. Let $\{x_0, x_1, \ldots, x_{\epsilon_1}, y_1, \ldots , y_{\epsilon_1}\}$ and $\{u_0, u_1, \ldots, u_{\epsilon_2}, w_1, \ldots , w_{\epsilon_2}\}$
be bases of $\fu_1$ and $\fu_2$, respectively. We define 
$J_\fu$ on $\fu_1\oplus \fu_2$ as the unique endomorphism satisfying:
\[ J_\fu x_0=u_0+w_1+x_1+y_2, \qquad J_\fu x_j=x_{j+1} +y_j+y_{j+2}, \qquad J_\fu u_k=u_{k+1}+w_k+w_{k+2},\]
\[
J_\fu u_0=-x_0-y_1+u_1+w_2, \qquad J_\fu y_j=-y_{j+1} -x_j, \qquad J_\fu w_k=-w_{k+1}-u_k
\]
for every $1\leq j\leq\epsilon_1,\ 1\leq k\leq\epsilon_2$ (with the convention that $x_j=y_j=0$ for $j>\epsilon_1$ and $u_k=w_k=0$ for $k>\epsilon_2$). It is easy to check that $J_\fu^2=-\Id_\fu$. Given two square matrices $M$ and $N$,  we denote by $M\oplus N$  the block matrix $\begin{pmatrix}
    M& \\ & N
\end{pmatrix}$. Set
\[J_\v := J_{\v_1}\oplus \cdots \oplus J_{\v_r}\oplus J_{\fu}, \quad \text{ where } \, J_{\v_\ell} =\begin{bmatrix}
    0&-\Id_\ell \\
    \Id_\ell &0
\end{bmatrix}.\] Here $\Id_\ell$ is the identity $k_\ell \times k_\ell$ or $2k_\ell \times 2k_\ell$ matrix, depending on whether $\alpha_\ell\in \R \cup \{\infty \}$ or $\alpha_\ell\in \C\setminus \R$.

Let $J_\z$ be  the endomorphism of $\z$ satisfying $J_\z z_1=-z_2$, $J_\z z_2=z_1$, and set $J:= J_\z \oplus J_\v$. 
 Clearly, $J^2= - \Id$ and we show next that $N_J=0$. 

The subspaces $\z,\v_1,\ldots,\v_r$ and $\fu$ are preserved by $J$ and mutually commute. Thus $N_J(x,y)=0$ whenever $x,y$ belong to two different subspaces in this list, or when $x,y\in\z$. Moreover, both $j(z_1)|_{\v_1 \oplus \cdots \oplus \v_r}$ and $  j(z_2)|_{\v_1 \oplus \cdots \oplus \v_r}$, commute with $J_{\v_1}\oplus \cdots \oplus J_{\v_r}$, which shows that $[x,y]=[Jx,Jy]$ and $[x,Jy]+[Jx,y]=0$, (and hence $N_J(x,y)=0$), whenever $x,y\in \z\oplus\v_1 \oplus \cdots \oplus \v_r$. 

By bilinearity, it remains to check that $N_J(x,y)$ vanishes whenever $x,y\in \fu$. Again, since the subspaces 
$\fu_{2i-1}\oplus\fu_{2i}$ are $J_\fu$-invariant for $1\leq i \leq m$ and mutually commute by the definition of the $j(z)$ maps, it is enough to check the vanishing of $N_J(x,y)$ for $x,y\in \fu_1\oplus\fu_2$. 

Finally, the general properties of the Nijenhuis tensor \eqref{nij} show that it suffices to check the vanishing of the following expressions:
$$N_J(x_0,y_j),\qquad N_J(x_0,w_k),\qquad N_J(y_i,y_j),\qquad N_J(y_i,w_k),\qquad N_J(w_k,w_l),$$
for $1\leq i,j\leq\epsilon_1,\ 1\leq k,l\leq\epsilon_2$.

From \eqref{eq:jz-odd-v}, the only non-zero brackets in $\fu_1\oplus\fu_2$ are 
$$[y_j,x_{j}]=z_1,\qquad [y_j,x_{j-1}]=z_2,\qquad [w_k,u_{k}]=z_1,\qquad [w_k,u_{k-1}]=z_2,$$
for $1\leq j\leq\epsilon_1,\ 1\leq k\leq\epsilon_2$. We compute for $j\geq1$:
\begin{eqnarray*}
    N_J(x_0,y_j)&=&[x_0,y_j]-[u_0+w_1+x_1+y_2,-y_{j+1} -x_j]\\&&+J([u_0+w_1+x_1+y_2,y_j]+[x_0,-y_{j+1} -x_j])\\
    &=&\delta_{1,j}(-z_2-Jz_1)+\delta_{2,j}(z_1-Jz_2)=0,
\end{eqnarray*}
and similarly for $k\geq 1$:
\begin{eqnarray*}N_J(x_0,w_k)&=&[x_0,w_k]-[u_0+w_1+x_1+y_2,-w_{k+1}-u_k]\\&&+J([u_0+w_1+x_1+y_2,w_k]+[x_0,-w_{k+1}-u_k])\\
&=&\delta_{1,k}(z_1-Jz_2)=0.\end{eqnarray*}
Then, for $1\leq i<j\leq\epsilon_1$:
\begin{eqnarray*}
    N_J(y_i,y_j)&=&[y_i,y_j]-[y_{i+1}+x_i,y_{j+1}+x_j]+J([-y_{i+1}-x_i,y_j]+[y_i,-y_{j+1}-x_j])\\
    &=&\delta_{i,{j-1}}(-z_1+Jz_2)=0,
\end{eqnarray*}
and for $1\leq k<l\leq\epsilon_2$:
\begin{eqnarray*}
    N_J(w_k,w_l)&=&[w_k,w_l]-[w_{k+1}+u_k,w_{l+1}+u_l]+J([-w_{k+1}-u_k,w_l]+[w_k,-w_{l+1}-u_l])\\
    &=&\delta_{{k},l-1}(-z_1+Jz_2)=0.
\end{eqnarray*}
Finally, for $1\leq i\leq\epsilon_1,\ 1\leq k\leq\epsilon_2$, we get directly:
\begin{eqnarray*}
    N_J(y_i,w_k)&=&[y_i,w_k]-[y_{i+1}+x_i,w_{k+1}+u_k]+J([-y_{i+1}-x_i,w_k]+[y_i,-w_{k+1}-u_k])=0.
\end{eqnarray*}

We note that in case $\mathcal S$ is as in \eqref{eq:S-abelian}, then  both $j(z_1)$ and $  j(z_2)$, commute with $J_\v$.  Then it follows from \cite[Lemma 5.6]{Bar} that $J:= J_\z \oplus J_\v$ is an abelian complex structure on $\n_\mathcal S$. 

It was shown in \cite[Section 3]{Lau-Osc} that the $2$-step nilpotent Lie algebras with $2$-dimensional center which are non-singular are those of the form  $\n_\mathcal S$ with $\mathcal S$ given by 
    \[  \mathcal S=\{(\alpha_1,k_1), \ldots ,(\alpha_r,k_r)  \}, \quad \alpha_\ell\in \C\setminus \R, \; k_\ell \in \N ,            \] 
 which is a particular case of \eqref{eq:S-abelian}.    Therefore, the last assertion in the statement follows.
    \end{proof}

As a direct consequence of this result we obtain Theorem \ref{thmC} stated in the introduction.

Let us now have a closer look at the $8$-dimensional $2$-step nilpotent Lie algebras which admit abelian complex structures.
As we mentioned in Section \ref{sec:prelim}, the dimension of the center of a $2$-step nilpotent Lie algebra $\n$ has to be even in order for $\n$ to admit abelian complex structures.

A further obstruction for $\n$ to admit abelian complex structures can be formulated in terms of its pencil $P_\n$.

\begin{lem}\label{lem:obstruction-abelian}
Let $\n$ be a $2$-step nilpotent Lie algebra. If the set $S(P_\n)$ of decomposable $2$-forms in $P_\n$ contains a plane $P$, then $\n$ admits  no abelian complex structure.
\end{lem}

\begin{proof}
Assume that $\n$ admits an abelian complex structure $J$ and let $\omega_1$ and $\omega_2$ be two linearly independent elements in $P$. Then $\omega_1+\omega_2\in P\subset S(P_\n)$, so 
$$\omega_1\wedge\omega_2=\frac12\left((\omega_1+\omega_2)\wedge(\omega_1+\omega_2)-\omega_1\wedge\omega_1-\omega_2\wedge\omega_2\right)=0.$$ We thus may write $\omega_i=\eta\wedge\zeta_i$, $i=1,2$ for some 1-forms $\eta,\zeta_1,\zeta_2\in\n^*$. Then $\ker\omega_i=\ker\eta\cap\ker\zeta_i$ are two distinct $J$-invariant hyperplanes of $\ker\eta$, so $\ker\eta=\ker\omega_1+\ker\omega_2$ is $J$-invariant. This is impossible, since $\dim\ker\eta=\dim\v-1$ is odd.
\end{proof}

Applying our results in dimension $8$ we then obtain the following classification.

\begin{thm}\label{thm:abelian-dim8}
An $8$-dimensional $2$-step nilpotent Lie algebra $\n$ admits an abelian complex structure if and only if $\dim\z=2$ or $\n$ is isomorphic to one of
\[
\R^5\oplus\h_3,\quad \R^3\oplus\h_5,\quad \R^2\oplus\h_3\oplus\h_3,\quad \R^2\oplus\h_3^\C,\quad \R^2\oplus\n_{6,2},
\]
\[
\R\oplus\n_{7,3,B_1},\quad \R\oplus\n_{7,3,D},\quad \R\oplus\n_{7,3,D_1},\quad \n_{8,4}^{1\text{-bis}}.
\]
\end{thm}

\begin{proof}
The case $\dim\z=2$ follows from Theorem \ref{thm:dimz=2}. For $\R\oplus\n_{7,3,B_1}$, $\R\oplus\n_{7,3,D}$, $\R\oplus\n_{7,3,D_1}$ and for $\n_{8,4}^{1\text{-bis}}$, the complex structures given in Tables~\ref{tab:complex-3} and~\ref{tab:complex-4} are abelian. The remaining Lie algebras in the statement are of the form $\R^2\oplus\mathfrak m$, where $\mathfrak m$ is a $6$-dimensional $2$-step nilpotent Lie algebra  admitting an abelian complex structure $J_{\mathfrak m}$ (see, for instance, \cite{ABD}). The direct sum of $J_{\mathfrak m}$ with the standard complex structure on $\R^2$ is then an abelian complex structure on $\R^2\oplus\mathfrak m$.

It remains to show that no other Lie algebra admits an abelian complex structure. The Lie algebras $\R^3\oplus\n_{5,2}$, $\R^2\oplus\n_{6,3}$, $\R\oplus\n_{7,2,A}$, $\R\oplus\n_{7,2,B}$, $\n_{5,2}\oplus\h_3$ and $\n_{8,3}^k$ have odd-dimensional center, so they admit no abelian complex structure. For the others, Lemma~\ref{lem:obstruction-abelian} applies, with the planes $P=\Span(e^{12},e^{23})$ for $\R\oplus\n_{7,3,A}$ and $\R\oplus\n_{7,3,B}$, the plane $P=\Span(e^{23},e^{24})$ for $\R\oplus\n_{7,3,C}$, and  the plane $P=\Span(e^{12},e^{13})$ for $\n_{8,4}^{1}$, $\n_{8,4}^{2}$ and $\n_{8,4}^{3}$.
\end{proof}

\medskip
\section{Appendix: The classification of $8$-dimensional $2$-step nilpotent Lie algebras}\label{app}

In this appendix, we summarize the classification of real $8$-dimensional $2$-step nilpotent Lie algebras $\n$, organized by $\dim(\n')$. Recall that $\n'\subseteq\mathfrak z(\n)$, and $\n$ has an abelian factor if and only if $\n'\subsetneq\mathfrak z(\n)$. The algebras with $\mathfrak z(\n)=\n'$ correspond to the algebras classified in \cite[Tables 1, 2 and 3]{BDG}. These include all the indecomposable ones, together with the decomposable algebras whose summands have no abelian factors. Since there is no $2$-step nilpotent Lie algebra without abelian factors of dimension $4$, and the only ones of dimension at most $5$ are $\mathfrak h_3$, $\mathfrak h_5$ and $\mathfrak n_{5,2}$, a dimension count shows that the decomposable ones of \cite{BDG} are exactly $\mathfrak h_3\oplus\mathfrak h_5$ and $\mathfrak n_{5,2}\oplus\mathfrak h_3$. The remaining algebras are of the form $\mathbb R\oplus\mathfrak m$, where $\mathfrak m$ is a $7$-dimensional $2$-step nilpotent Lie algebra, and can be obtained from the classification in \cite{Gong}.

For each Lie algebra $\n$, the structure equations are written with respect to a basis $\{e^1,\ldots,e^8\}$ of the dual space $\n^*$. As usual, $e^{ij}$ stands for $e^i\wedge e^j$.

The notation $\mathfrak n_{n,n'}$ or $\mathfrak n_{n,n',\bullet}$ means that the Lie algebra has dimension $n$ and derived algebra of dimension $n'$, while different capital letters in the third argument are used to distinguish non-isomorphic Lie algebras whose derived algebras have the same dimension. The indecomposable $8$-dimensional algebras are denoted by $\mathfrak n_{8,n'}^{k}$, where $k$ is the number of the corresponding orbit in \cite[Tables 1, 2 and 3]{BDG}, and $\mathfrak n_{8,n'}^{k\text{-bis}}$ denotes the second real form within the same complex orbit. Note that algebras $\n_{8,2}^5$ in \cite[Table 1]{BDG} and $\n^{11}_{8,3}$ in \cite[Table 2]{BDG} do not appear in the tables below since they are isomorphic to the decomposable algebras $\mathfrak h_3\oplus\mathfrak h_5$ and $\mathfrak n_{5,2}\oplus\mathfrak h_3$, respectively.
Moreover, $\mathfrak h_n$ denotes the Heisenberg Lie algebra of dimension $n$ (for odd $n$), and $\mathfrak h_3^{\mathbb C}$ denotes the real Lie algebra underlying the complex $3$-dimensional Heisenberg Lie algebra.

\begin{itemize}
\item $8$-dimensional $2$-step nilpotent Lie algebras $\mathfrak n$ with $\dim(\mathfrak n')=1$:
\begin{align*}
\Lb{\mathbb R^5 \oplus\mathfrak h_3}&= \Rb{\big(0,0,0,0,0,0,0,e^{12}\big),}\\
\Lb{\mathbb R^3 \oplus\mathfrak h_5}&= \Rb{\big(0,0,0,0,0,0,0,e^{12}+e^{34}\big),}\\
\Lb{\mathbb R \oplus\mathfrak h_7}&= \Rb{\big(0,0,0,0,0,0,0,e^{12}+e^{34}+e^{56}\big).}
\end{align*}
\end{itemize}
\begin{itemize}
\item $8$-dimensional $2$-step nilpotent Lie algebras $\mathfrak n$ with $\dim(\mathfrak n')=2$:
\begin{align*}
\Lb{\mathbb R^3 \oplus\mathfrak n_{5,2}} &= \Rb{\big(0,0,0,0,0,0,e^{12},e^{13}\big),}\\
\Lb{\mathbb R^2 \oplus\mathfrak h_3\oplus\mathfrak h_3}&= \Rb{\big(0,0,0,0,0,0,e^{12},e^{34}\big),}\\
\Lb{\mathbb R^2 \oplus\mathfrak h_3^{\mathbb C}}&= \Rb{\big(0,0,0,0,0,0,e^{13}-e^{24},e^{14}+e^{23}\big),}\\
\Lb{\mathbb R^2 \oplus\mathfrak n_{6,2}} &= \Rb{\big(0,0,0,0,0,0,e^{12},e^{14}+e^{23}\big),}\\
\Lb{\mathbb R \oplus\mathfrak n_{7,2,A}} &= \Rb{\big(0,0,0,0,0,0,e^{12},e^{14}+e^{35}\big),}\\
\Lb{\mathbb R \oplus\mathfrak n_{7,2,B}} &= \Rb{\big(0,0,0,0,0,0,e^{12}+e^{34},e^{15}+e^{23}\big),}\\
\Lb{\mathfrak h_3\oplus\mathfrak h_5} &= \Rb{\big(0,0,0,0,0,0,e^{12},e^{34}+e^{56}\big),}\\
\Lb{\mathfrak n_{8,2}^{1}} &= \Rb{\big(0,0,0,0,0,0,e^{12}+e^{56},e^{34}+e^{56}\big),}\\
\Lb{\mathfrak n_{8,2}^{1\text{-bis}}} &= \Rb{\big(0,0,0,0,0,0,e^{12}-e^{36}+e^{45},e^{34}-e^{56}\big),}\\
\Lb{\mathfrak n_{8,2}^{2}} &= \Rb{\big(0,0,0,0,0,0,e^{15}+e^{46},e^{14}+e^{23}\big),}\\
\Lb{\mathfrak n_{8,2}^{3}} &= \Rb{\big(0,0,0,0,0,0,e^{15}+e^{26}+e^{34},e^{14}+e^{23}\big),}\\
\Lb{\mathfrak n_{8,2}^{4}} &= \Rb{\big(0,0,0,0,0,0,e^{13}+e^{25},e^{16}+e^{24}\big),}\\
\Lb{\mathfrak n_{8,2}^{6}} &= \Rb{\big(0,0,0,0,0,0,e^{16}+e^{25}+e^{34},e^{12}\big).}
\end{align*}
\end{itemize}
\begin{itemize}
\item $8$-dimensional $2$-step nilpotent Lie algebras $\mathfrak n$ with $\dim(\mathfrak n')=3$:
{\allowdisplaybreaks \begin{align*}
\Lb{\mathbb R^2 \oplus\mathfrak n_{6,3}} &= \Rb{\big(0,0,0,0,0,e^{12},e^{13},e^{23}\big),}\\
\Lb{\mathbb R \oplus\mathfrak n_{7,3,A}}&= \Rb{\big(0,0,0,0,0,e^{12},e^{23},e^{24}\big),}\\
\Lb{\mathbb R \oplus\mathfrak n_{7,3,B}}&= \Rb{\big(0,0,0,0,0,e^{12},e^{23},e^{34}\big),}\\
\Lb{\mathbb R \oplus\mathfrak n_{7,3,B_1}}&= \Rb{\big(0,0,0,0,0,e^{12}-e^{34},e^{13}+e^{24},e^{14}\big),}\\
\Lb{\mathbb R \oplus\mathfrak n_{7,3,C}}&= \Rb{\big(0,0,0,0,0,e^{12}+e^{34},e^{23},e^{24}\big),}\\
\Lb{\mathbb R \oplus\mathfrak n_{7,3,D}} &= \Rb{\big(0,0,0,0,0,e^{12}+e^{34},e^{13},e^{24}\big),}\\
\Lb{\mathbb R \oplus\mathfrak n_{7,3,D_1}} &= \Rb{\big(0,0,0,0,0,e^{12}-e^{34},e^{13}+e^{24},e^{14}-e^{23}\big),}\\
\Lb{\mathfrak n_{5,2}\oplus\mathfrak h_3} &= \Rb{\big(0,0,0,0,0,e^{12},e^{13},e^{45}\big),}\\
\Lb{\mathfrak n_{8,3}^{1}} &= \Rb{\big(0,0,0,0,0,e^{12}+e^{45},e^{25}+e^{34},e^{15}+e^{23}\big),}\\
\Lb{\mathfrak n_{8,3}^{1\text{-bis}}} &= \Rb{\big(0,0,0,0,0,e^{12}+2e^{23}+2e^{45},e^{25}-2e^{34},-e^{15}+2e^{24}-2e^{35}\big),}\\
\Lb{\mathfrak n_{8,3}^{2}} &= \Rb{\big(0,0,0,0,0,e^{14}+e^{35},e^{15}+e^{24},e^{23}\big),}\\
\Lb{\mathfrak n_{8,3}^{2\text{-bis}}} &= \Rb{\big(0,0,0,0,0,e^{14}+e^{24}-e^{35},-e^{15}+e^{25}+e^{34},-e^{23}\big),}\\
\Lb{\mathfrak n_{8,3}^{3}} &= \Rb{\big(0,0,0,0,0,e^{15}+e^{34},e^{23},e^{14}+e^{24}\big),}\\
\Lb{\mathfrak n_{8,3}^{3\text{-bis}}} &= \Rb{\big(0,0,0,0,0,e^{15}+e^{34},e^{14}-e^{23},e^{13}+e^{24}\big),}\\
\Lb{\mathfrak n_{8,3}^{4}} &= \Rb{\big(0,0,0,0,0,e^{25}+e^{34},e^{15}+e^{23},e^{12}\big),}\\
\Lb{\mathfrak n_{8,3}^{5}} &= \Rb{\big(0,0,0,0,0,e^{13},e^{24},e^{15}+e^{25}\big),}\\
\Lb{\mathfrak n_{8,3}^{5\text{-bis}}} &= \Rb{\big(0,0,0,0,0,e^{13},e^{15}-e^{24},e^{14}+e^{25}\big),}\\
\Lb{\mathfrak n_{8,3}^{6}} &= \Rb{\big(0,0,0,0,0,e^{24},e^{15}+e^{23},e^{13}\big),}\\
\Lb{\mathfrak n_{8,3}^{7}} &= \Rb{\big(0,0,0,0,0,e^{15}+e^{24},e^{14}+e^{23},e^{13}\big),}\\
\Lb{\mathfrak n_{8,3}^{8}} &= \Rb{\big(0,0,0,0,0,e^{25}+e^{34},e^{12},e^{15}\big),}\\
\Lb{\mathfrak n_{8,3}^{9}} &= \Rb{\big(0,0,0,0,0,e^{13}+e^{25},e^{15}+e^{24},e^{12}\big),}\\
\Lb{\mathfrak n_{8,3}^{10}} &= \Rb{\big(0,0,0,0,0,e^{14}+e^{35},e^{23},e^{12}\big),}\\
\Lb{\mathfrak n_{8,3}^{12}} &= \Rb{\big(0,0,0,0,0,e^{15}+e^{23},e^{14},e^{12}\big).}
\end{align*}}
\end{itemize}
\begin{itemize}
\item $8$-dimensional $2$-step nilpotent Lie algebras $\mathfrak n$ with $\dim(\mathfrak n')=4$:
\begin{align*}
\Lb{\mathfrak n_{8,4}^{1}} &= \Rb{\big(0,0,0,0,e^{12},e^{34},e^{13},e^{24}\big),}\\
\Lb{\mathfrak n_{8,4}^{1\text{-bis}}} &= \Rb{\big(0,0,0,0,e^{12},e^{34},e^{13}-e^{24},-e^{14}-e^{23}\big),}\\
\Lb{\mathfrak n_{8,4}^{2}} &= \Rb{\big(0,0,0,0,e^{13},e^{24},e^{14}+e^{23},e^{12}\big),}\\
\Lb{\mathfrak n_{8,4}^{3}} &= \Rb{\big(0,0,0,0,e^{14},e^{23},e^{12},e^{13}\big).}
\end{align*}
\end{itemize}

\smallskip

\subsection*{Disclosure on the use of AI}
Claude by Anthropic was used by the third named author to make some preliminary numerical explorations about the existence of complex structures on some $2$-step nilpotent Lie algebras, to assist with the LaTeX formatting of the tables and to improve the Introduction. The rest of the typesetting, the proofs and the ideas behind them are all entirely produced by the authors.

\smallskip



\begin{thebibliography}{WWWW}


\bibitem{ABB}
A. Andrada, M. L. Barberis, B. Brienza, 
\emph{Holonomy of the Obata connection on $2$-step hypercomplex nilmanifolds}, Bull. Lond. Math. Soc. \textbf{58} (2026), Article ID e70338, 23 p.

\bibitem{ABD}
A. Andrada, M. L. Barberis, I. Dotti, 
\emph{Classification of abelian complex structures on $6$-dimensional Lie algebras},
J. Lond. Math. Soc., II. Ser. \textbf{83} (2011), 232--255.

\bibitem{AD} A. Andrada, I. Dotti, 
\emph{Killing-Yano $2$-forms on $2$-step nilpotent Lie groups}, Geom. Dedicata \textbf{212} (2021), 415--424.

\bibitem{AVe}
A. Andrada, S. Vera, \emph{
Complex structures on $2$-step nilpotent Lie algebras arising from graphs}, 
J. Pure Appl. Algebra \textbf{230} (2026), Paper No. 8, 39 p.

\bibitem{AV}
A. Andrada, R. Villacampa, \emph{Abelian balanced Hermitian structures on unimodular Lie algebras}, Transform. Groups \textbf{21} (2016), 903--927.

\bibitem{AN}
R. M. Arroyo, M. Nicolini, \emph{SKT structures on nilmanifolds}, Math. Z. \textbf{302} (2022), 1307--1320.

\bibitem{Bar}
M. L. Barberis, 
\emph{Complex structures on two-step nilpotent Lie groups}, preprint
arXiv:2508.05885 (2025).

\bibitem{BR}
G. Bazzoni, J. Rojo, 
\emph{$8$-dimensional $2$-step nilpotent Lie algebras over algebraically closed fields with $\operatorname{char}\neq2,3$}, Linear Algebra Appl. {\bf 752} (2027), 536--585.

\bibitem{BG}
C. Benson, C. Gordon, 
\emph{Kähler and symplectic structures on nilmanifolds}, Topology \textbf{27} (1988),
513--518.

\bibitem{BDG}
M. Borovoi, B. A. Dina, W. A.  de Graaf, 
\emph{Real non-degenerate two-step nilpotent Lie algebras of dimension eight}, 
Eur. J. Math. \textbf{10} (2024),  Paper No. 16, 27p. 

\bibitem{COUV}
M. Ceballos, A. Otal, L. Ugarte, R. Villacampa, \emph{Invariant complex structures on $6$-nilmanifolds: classification,
Frölicher spectral sequence and special Hermitian metrics}, J. Geom. Anal. {\bf 26} (2016), 252--286.

\bibitem{CdBM}
A. Clarke, V. del Barco, A. Moreno,
\emph{$\mathrm{G}_2$-instantons on $2$-step nilpotent Lie groups},
Adv. Theor. Math. Phys. {\bf 29} (2025), 751--813.

\bibitem{CFGU}
L. Cordero, M. Fernández, A. Gray,  L. Ugarte, \emph{Compact nilmanifolds with nilpotent complex structures: Dolbeault cohomology}, Trans. Amer. Math. Soc. \textbf{352} (2000), 5405--5433.

\bibitem{DR}
E. Damek, F. Ricci, \emph{A class of nonsymmetric harmonic Riemannian spaces}, Bull. Amer. Math. Soc., New Ser. {\bf 27} (1992), 139--142.

\bibitem{DM}
S. Dani, M. Mainkar, \emph{Anosov automorphisms on compact nilmanifolds associated with graphs},
Trans. Amer. Math. Soc. \textbf{357} (2005), 2235--2251. 

\bibitem{dBM}
 V. del Barco, A. Moroianu, \emph{Killing forms on $2$-step nilmanifolds}, J. Geom. Anal. \textbf{31} (2021), 863--887.

\bibitem{dBMR}
V. del Barco, A. Moroianu, A. Raffero,
\emph{Purely coclosed $\mathrm{G}_2$-structures on $2$-step nilpotent Lie groups}, Rev. Mat. Complut. {\bf 35} (2022), 323--359.

\bibitem{De}
J. Der\'e, \emph{Orthogonal bi-invariant complex structures on metric Lie
algebras}, {Ann. Global Anal. Geom. } \textbf{ 59} (2021), 157--177. 

\bibitem{DLV}
A. J. Di Scala, J. Lauret, L. Vezzoni,  \emph{Quasi-Kähler Chern-flat manifolds and complex $2$-step nilpotent
Lie algebras}, {Ann. Sc. Norm. Super. Pisa Cl. Sci. (5)} \textbf{11} (2012), 41--60.

\bibitem{DF}
 I. Dotti, A. Fino, \emph{Hypercomplex eight-dimensional nilpotent Lie groups}, J. Pure Appl. Algebra \textbf{184} (2003), 41--57.
 
\bibitem{Eb}
P. Eberlein, 
\emph{Geometry of $2$-step nilpotent groups with a left invariant metric}, 
Ann. Sci. Éc. Norm. Supér. (4) \textbf{27} (1994), 611--660.

\bibitem{Eb2}
P. Eberlein, \emph{Geometry of $2$-step nilpotent Lie groups with a left-invariant metric II}, Trans. Amer. Math. Soc.
{\bf 343} (1994), 805--828.

\bibitem{EFV1}
N. Enrietti, A. Fino, L. Vezzoni, \emph{Tamed symplectic forms and strong Kähler with torsion metrics},
J. Symplect. Geom. \textbf{10} (2012), 203--223.

\bibitem{EFV2}
N. Enrietti, A. Fino,  L. Vezzoni, \emph{The pluriclosed flow on nilmanifolds and tamed symplectic forms}, J.
Geom. Anal. \textbf{25} (2015), 883--909.

\bibitem{FGR}
A. Fino, G. Grantcharov, G. Russo,
\emph{On the $\nu$-invariant of two-step nilmanifolds with closed $\mathrm{G}_2$-structure}, preprint
arXiv:2409.06870 (2024).

\bibitem{FPS}
A. Fino, M. Parton, S. Salamon, \emph{Families of strong KT structures in six dimensions}, Comment. Math. Helv. \textbf{79} (2004), 317--340.

\bibitem{FV}
A. Fino, L. Vezzoni, \emph{On the existence of balanced and SKT metrics on nilmanifolds}, Proc.
Amer. Math. Soc. \textbf{144} (2016), 2455--2459.

\bibitem{Fol}
 G. Folland, 
 \emph{Subelliptic estimates and function spaces on nilpotent Lie groups}, Ark. Mat. {\bf 13} (1975), 161--207.

\bibitem{GZZ}
Q. Gao, Q. Zhao, F. Zheng, \emph{Maximal nilpotent complex structures}, Transform. Groups \textbf{28} (2023),
241--284.

\bibitem{Gau}
M. Gauger, \emph{On the classification of metabelian Lie algebras}, Trans. Amer. Math. Soc. \textbf{179} (1973), 293--329.

\bibitem{Gong}
M.-P. Gong, 
\emph{Classification of Nilpotent Lie Algebras of Dimension $7$ (Over Algebraically Closed Fields and $\mathbb R$).} PhD
thesis, University of Waterloo (Canada), 1998. 

\bibitem{Gor}
C. Gordon, \emph{Naturally reductive homogeneous Riemannian manifolds}, Canad. J. Math. {\bf 37} (1985), 467--487.

\bibitem{GR}
M. Goze, E. Remm, \emph{Non existence of complex structures on filiform Lie algebras}, 
Commun. Algebra {\bf 30} (2002), 3777--3788.

\bibitem{Kap}
A. Kaplan, \emph{Riemannian Nilmanifolds Attached to Clifford Modules}, Geom. Dedicata  \textbf{11} (1981), 
 127--136.

\bibitem{LU}
A. Latorre, L. Ugarte,
\emph{Nilmanifolds with non-nilpotent complex structures and their pseudo-Kähler geometry}, Rev. Real Acad. Cienc. Exactas Fis. Nat. Ser. A-Mat. {\bf 119} (2025), Paper No. 47, 29p.

\bibitem{LU2}
A. Latorre, L. Ugarte,
\emph{On astheno-K\"ahler nilmanifolds with balanced metrics}, preprint
arXiv:2606.28545 (2026).

\bibitem{LUV1}
A. Latorre, L. Ugarte, R. Villacampa, \emph{On the Bott-Chern cohomology and balanced Hermitian nilmanifolds}, Internat. J. Math. {\bf 25} (2014), Article ID 1450057. 

\bibitem{LUV2} A. Latorre, L. Ugarte,  R. Villacampa, \emph{The ascending central series of nilpotent Lie algebras with complex structure}, Trans. Amer. Math. Soc. {\bf 372} (2019), 3867--3903. 

\bibitem{LUV3}
A. Latorre, L. Ugarte, R. Villacampa,
\emph{Complex structures on nilpotent Lie algebras with one-dimensional center},
J. Algebra \textbf{614} (2023), 271--306.

\bibitem{Lauret}
J. Lauret, \emph{Homogeneous nilmanifolds attached to representations of compact Lie groups}, Manuscripta Math.
{\bf 99} (1999), 287--309.

\bibitem{Lau-Osc}
 J. Lauret, D. Oscari, 
\emph{On nonsingular $2$-step nilpotent Lie algebras}, 
Math. Res. Lett. \textbf{21} (2014),  553--583. 

\bibitem{LT}
F. Levstein, A. Tiraboschi, \emph{Classes of $2$-step nilpotent Lie algebras}, Comm. Alg. \textbf{27} (1999), 2425--2440.

\bibitem{LoG}
E. Lo Giudice, \emph{Positive Hermitian curvature flow on $2$-step nilpotent Lie groups},
Ann. Global Anal. Geom. {\bf 69} (2026), Paper No. 6, 21 p.

\bibitem{LRV}
B. Luporini, S. Reggiani, F. Vittone,
\emph{On the structure of $2$-step nilpotent Lorentzian naturally reductive Lie groups}, preprint
arXiv:2607.18017 (2026).

\bibitem{MO}
A. Moroianu, G.P. Ovando, 
\emph{Magnetic Killing tensors and first integrals of the magnetic flow}, Lett. Math. Phys. {\bf 116} (2026), Paper No. 65, 26 p.

\bibitem{MRV}
A. Moroianu, A. Raffero, L. Vezzoni,
\emph{The heterotic $\mathrm{G}_2$-system on $2$-step nilmanifolds endowed with principal torus bundles}, preprint
arXiv:2512.16817 (2025), to appear in Adv. Theor. Math. Phys.

\bibitem{OS}
G.P. Ovando, M. Subils, \emph{Symplectic structures on low dimensional $2$-step nilmanifolds}, Mat. Contemp. {\bf 60} (2024), 225--254.

\bibitem{Par}
S. Pareja, \emph{Geometría casi-compleja}, Bachelor's thesis, Universidad de Zaragoza (2022). Available at \texttt{https://zaguan.unizar.es/record/124830}.

\bibitem{PV}
L. Pizarro, R. Villacampa,
\emph{On abelian complex structures on nilpotent Lie algebras}, preprint
arXiv:2608.04914 (2026).

\bibitem{Rol}
S. Rollenske, \emph{Geometry of nilmanifolds with left-invariant complex structure and deformations in the large}, {Proc.
London Math. Soc.} \textbf{99} (2009),  425--460.

\bibitem{Sal}
S. M. Salamon, 
\emph{Complex structures on nilpotent Lie algebras}, J. Pure Appl. Algebra
\textbf{157} (2001), 311--333. 

\bibitem{ST}
T. Sferruzza, N. Tardini, \emph{$p$-Kähler and balanced structures on nilmanifolds with nilpotent complex
structures}, Ann. Global Anal. Geom. \textbf{62} (2022), 869--881.

\bibitem{Snow}
J. E. Snow, \emph{ Invariant complex structures on four dimensional solvable real Lie groups},
Manuscripta Math. \textbf{66} (1990), 397--412.

\bibitem{UV}
L. Ugarte, R. Villacampa, \emph{Balanced Hermitian geometry on $6$-dimensional nilmanifolds}, Forum Math.
\textbf{27} (2015), 1025--1070. 

\bibitem{Zhang}
J. Zhang, \emph{Complex structures on stratified Lie algebras}, Bull. Aust. Math. Soc. {\bf 106} (2022), 320--332.

\end{thebibliography}
\end{document}